\documentclass{amsart}

\usepackage{amssymb}
\usepackage{supertabular}
\usepackage{url}

\theoremstyle{plain}
\newtheorem{thm}{Theorem}[section]
\newtheorem{cor}[thm]{Corollary}
\newtheorem{lemma}[thm]{Lemma}
\newtheorem{claim}{Claim}

\theoremstyle{definition}
\newtheorem{defn}[thm]{Definition}
\newtheorem{eg}[thm]{Example}
\newtheorem{rmk}[thm]{Remark}

\numberwithin{equation}{section}
\numberwithin{table}{section}

\def\Frak{\mathfrak}
\def\Cal{\mathcal}
\def\Bb{\mathbb}

\def\A{{\Bb A}}
\def\F{{\Bb F}}
\def\bP{{\Bb P}}
\def\Q{{\Bb Q}}
\def\bS{{\Bb S}}
\def\Z{{\Bb Z}}
\def\O{{\Cal O}}
\def\H{{\Bb H}}

\def\C{{\Bb C}}

\def\q{{\lambda}}
\def\p{{\pi}}
\def\k{{\kappa}}
\def\i{{i}}

\def\fR{{\Frak R}}

\def\qtimes{\otimes_{\Q}}
\def\dtimes{\otimes_{D}}
\def\Ql{{\Q_\ell}}
\def\orho{{\overline\rho}}
\def\rhox{{\rho_{2, X}}}
\def\orhox{{\overline\rho_{2, X}}}
\def\orhoy{{\overline\rho_{2, Y}}}
\def\divides{{\mid}}

\def\mult{^\times}
\def\inv{^{-1}}

\def\set#1{\left\{#1\right\}}
\def\sset#1{\{#1\}}

\DeclareMathOperator{\End}{End}
\DeclareMathOperator{\Aut}{Aut}
\DeclareMathOperator{\Mat}{Mat}
\DeclareMathOperator{\Gal}{Gal}
\DeclareMathOperator{\GL}{GL}
\DeclareMathOperator{\SL}{SL}
\DeclareMathOperator{\Sp}{Sp}
\DeclareMathOperator{\PGL}{PGL}
\DeclareMathOperator{\PSL}{PSL}
\DeclareMathOperator{\discr}{\delta}
\DeclareMathOperator{\chr}{char}
\DeclareMathOperator{\Hom}{Hom}
\DeclareMathOperator{\Perm}{Perm}
\DeclareMathOperator{\rank}{rank}
\DeclareMathOperator{\res}{res}
\DeclareMathOperator{\Id}{Id}

\begin{document}

\title{Hyperelliptic jacobians with real multiplication}

\author{Arsen Elkin}
%\email{arsen\char`\@psu.edu}

\address{Department of Mathematics, The Pennsylvania State University,
University Park, PA 16802, USA}

\email{arsen@psu.edu}



\date{\today}

\begin{abstract}
Let $K$ be a field of characteristic $p \neq 2$, 
and let $f(x)$ be a sextic polynomial irreducible over
$K$ with no repeated roots, whose Galois group is isomorphic to $\A_5$.
If the jacobian $J(C)$ of the hyperelliptic curve $C:y^2=f(x)$ admits 
real multiplication over the ground field from an order of a real
quadratic field $D$, then either its endomorphism algebra is isomorphic 
to $D$, or $p > 0$ and $J(C)$ is a supersingular abelian variety. 
The supersingular outcome cannot occur when $p$ splits in $D$.
\end{abstract}

\maketitle
%===================================================================
\section{Statement of Results}
%===================================================================
Let $K$ be a field and let $K_a$ be its fixed algebraic closure.
Let $f(x)\in K[x]$ be an irreducible polynomial
of degree $n=6$ with no repeated roots.
Denote by $\fR_f\subset K_a$ the set of roots of $f$ and 
by $K(\fR_f)$ the extension of $K$ generated by the 
elements of $\fR_f$, that is, the splitting  field of $f$ over $K$.
We write $\Gal(f/K)$ for the Galois group of $K(\fR_f)/K$, or simply $\Gal(f)$ when
no confusion over the ground field arises. This
group acts on the elements of $\fR_f$ by permutations, 
and it is well-known that this action is transitive if and only if
$f$ is irreducible over $K$.

Consider the hyperelliptic curve defined over $K$ by
$$
C_f: y^2=f(x).
$$
Let $J(C_f)$ be its jacobian, 
$\End(J(C_f))$ the ring of $K_a$-endo\-morphisms of $J(C_f)$, 
and $\End_K(J(C_f))$ the ring of $K$-endomorphisms 
of $J(C_f)$. 
We define algebras $\End^0(J(C_f)):=\End(J(C_f))\otimes\Q$
and $\End_K^0(J(C_f)):=\End_K(J(C_f))\otimes\Q$.
One may view $\End(J(C_f))$ and $\End_K(J(C_f))$ 
as orders in the corresponding $\Q$-algebras.

For every algebraic curve $C$, the ring $\End(J(C))$ contains 
the multiplications by integers; that is, $\Z\cdot\Id_{J(C)}\subset\End(J(C))$,
where $\Id_{J(C)}$ is the identity automorphism of $J(C)$.
Examples with the precise equality $\End(J(C))=\Z$ are 
harder to find (see \cite{MasserSHJ, MoriER1, MoriER2}). 
In \cite{ZarCM, ZarCMP, ZarNS} Yu.~Zarhin 
proves that if $\Gal(f)$ is isomorphic to either
$\A_n$ or $\bS_n$, then this equality holds for 
the curve $C_f$, or $\chr(K)=3$, $n=5$ or $6$, and $J(C_f)$ is 
a supersingular abelian variety.
In \cite{ZarVS2R} he proves that if the roots $\fR_f$ of a polynomial $f$ 
can be identified with $\bP^{m-1}(\F_q)$ 
for some odd prime power $q$ and integer $m>2$,
and $\Gal(f)$ contains $\PSL_m(\F_q)$ as a subgroup, then
$\End(J(C_f))=\Z$ or $J(C_f)$ is a supersingular abelian variety.
This statement is not necessarily true for $m=2$.

Similarly, examples of hyperelliptic curves
$C_f$ with $\End(J(C_f))$ containing an order of
a real quadratic field (admitting real multiplication)
are known (see \cite{WilEM, WilRM}).
The purpose of this paper is to provide a tool for construction of
explicit examples of hyperelliptic curves $C_f$ for which $\End(J(C_f))$ is
isomorphic to the ring of integers of a real quadratic field $D=\Q(\omega)$,
where $\omega^2=d$ for some square-free integer $d \geq 2$
called the {\em reduced discriminant} of $D$. Simply put, $\omega=\sqrt{d}$.
Consider polynomials $f$ of degree $n=6$,
in which case curves $C_f$ have genus $g=2$. 
We will denote the discriminant of an algebraic number field $D$ 
(or its order $\O$) by $\discr(D)$ (respectively $\discr(\O)$) .
We fix all of the above notation throughout the paper. 

We will prove the following statement.
\begin{thm}[The Main Theorem]\label{mainthm}
Let $K$ be a field of characteristic $p \neq 2$.
Let $f(x)\in K[x]$ be an irreducible separable polynomial of degree $n=6$.
Let $J(C_f)$ be the jacobian of the hyperelliptic curve $C_f: y^2=f(x)$.
Let $D$ be a real quadratic field.
Assume that $\Gal(f)$ and $J(C_f)$ enjoy the following properties:
\renewcommand{\theenumi}{\alph{enumi}}
\begin{enumerate}
\item $\Gal(f)\cong\A_5$,
\item there exists an injective ring homomorphism 
$\i: D \hookrightarrow \End_K^0(J(C_f))$ such that
$\i(1)=\Id_{J(C_f)}$, the identity automorphism of $J(C_f)$.
\end{enumerate}
Then $\End_K(J(C_f))$ is isomorphic to an order of $D$ with
$$
\discr(\End_K(J(C_f)))\equiv5\bmod8,
$$
and one of the following conditions holds:
\renewcommand{\theenumi}{\roman{enumi}}
\begin{enumerate}
\item\label{nonsup}
$J(C_f)$ is an absolutely simple abelian variety and 
$$
\End(J(C_f))=\End_K(J(C_f)).
$$
\item\label{supsing}
$p > 0$ and $J(C_f)$ is a supersingular abelian variety. 
Moreover, $p$ does not split in $D$.
\end{enumerate}
\end{thm}

%----------------------------------------------------------------

This theorem is a modification of the results in \cite{ZarCM, ZarCMP, ZarNS}.
The structure of the paper is as follows.
The proof of Theorem \ref{mainthm} will be given in Section \ref{mainpf}. 
In Section \ref{NSSect}, the impossibility of the 
supersingular outcome when $p=\chr(K)$ splits in $D$
is proven. Examples in characteristic zero will be given in Section \ref{examples}.
Examples, both supersingular and not, in positive characteristic will be given
in Section \ref{PosEg}.

\begin{rmk}
If $\O$ is an order of $D = \Q(\omega)$ with conductor $c$ (see Section \ref{mainpf}), 
then the condition $\discr(\O)\equiv5\bmod8$ is equivalent
to $d\equiv5\bmod8$ and $c$ being odd.
\end{rmk}

\begin{rmk}
If $D$ is a quadratic field with reduced discriminant $d$, then in order for an odd 
rational prime $p$ not to split in $D$ it is necessary and sufficient that either
$p\divides d$ (in which case, $p$ ramifies in $D$) or $\left(\frac{d}{p}\right)=-1$ 
(that is, $p$ is inert in $D$).
\end{rmk}

\begin{rmk}\label{dtrmk}
Since $f$ is irreducible over $K$, the action of $\Gal(f)$ on the 
set $\fR_f$ on the roots of $f$ is transitive. 
According to \cite[Table. 2.1, p. 60]{DixPG},
the six roots of an irreducible polynomial $f$ with $\Gal(f)\cong\A_5$ 
can be identified with the elements of $\bP^1(\F_5)$ in such a 
way that the action of $\Gal(f)$ on the set $\fR_f$ is isomorphic to 
the action of $\PSL_2(\F_5)\cong\A_5$
on $\bP^1(\F_5)$. 
Note that this action is doubly transitive.
\end{rmk}
%------------------------------------------------------------

In proving Theorem \ref{mainthm} we will use the following statement,
whose proof will be given in Section \ref{auxpf}. 
\begin{thm}\label{auxthm}
Let $K$ be a field of characteristic different from $2$.
Let $f(x)\in K[x]$ be an irreducible separable polynomials
of degree $n=6$ with $\Gal(f)\cong\A_5$.
Let $X=J(C_f)$ be the jacobian of the hyperelliptic curve $C_f: y^2=f(x)$.
Let $R$ be a subalgebra of $\End_{\F_2}(X_2)$ containing the
identity automorphism of $X_2$ such that
$$
{^\sigma}\!u \in R\quad\text{for each }\,u \in R, \sigma\in\Gal(K), 
$$
where
$$
{^\sigma}\!u: x\mapsto \sigma u(\sigma\inv x),\quad x\in X_2.
$$
Then $\End_{\Gal(K)}(X_2) \cong \F_4$ and
one of the following $\F_2$-algebra isomorphisms holds:
\renewcommand{\theenumi}{\roman{enumi}}
\begin{enumerate}
\item\label{csf2} $R \cong \F_2$.
\item\label{csf4} $R \cong \F_4$.
\item\label{csm2f4} $R \cong \Mat_2(\F_4)$.
\item\label{csm4f2} $R = \End_{\F_2}(X_2)\cong \Mat_4(\F_2)$.
\end{enumerate}
\end{thm}
%----------------------------------------------
As a corollary to this theorem we will also obtain
\begin{thm}\label{rankthm}
Let $K$ be a field of characteristic different from $2$.
Let $f(x)\in K[x]$ be an irreducible separable polynomial 
of degree $n=6$ with $\Gal(f)\cong\A_5$. 
Let $C_f:y^2=f(x)$ be a hyperelliptic curve over $K$.
Then either
\begin{itemize}
\item[(a)] $\End_K(J(C_f)) = \Z$, or
\item[(b)] $\End_K(J(C_f))$ is isomorphic to an order of a quadratic field with
$$
\discr(\End_K(J(C_f)))\equiv5\bmod8.
$$
\end{itemize}
%In particular, the rank of $\End_K(J(C))$ over $\Z$ does not exceed $2$.
\end{thm}
The proof of this theorem will be given in Section \ref{mainpf}.

%----------------------------------------------
The following statement will be used in order to show that
jacobians we produce as examples in Section \ref{examples} are pairwise non-isogenous. 
Its proof will be given in Section \ref{nonisogenuity}.
\begin{thm}\label{NICor}
Let $K$ be a field of characteristic $0$, $f(x)$, $h(x)\in K(x)$ 
irreducible separable polynomials over $K$ of degree $n=6$ each, 
such that $K(\fR_f)$ and $K(\fR_h)$ are linearly disjoint 
extensions of $K$. Assume that $J(C_f)$ and $J(C_h)$
satisfy the conditions of Theorem \ref{mainthm}.
Then 
$$
\Hom(J(C_f), J(C_h))=0
\quad\text{and}\quad
\Hom(J(C_h), J(C_f))=0.
$$
\end{thm}

%=======================================================
\section{Proof of the Main Result}
\label{mainpf}
%=======================================================
For an abelian variety $X$ defined over 
a field of characteristic distinct from $2$,
the natural action of $\End(X)$ on $X_2$ induces an injective homomorphism
$$
\End(X) \otimes \Z/2\Z \hookrightarrow \End(X_2).
$$
By \cite[p. 501]{SerTateRed} we have $\End_K(X) = \End_K^0(X) \cap \End(X)$, 
so the map
$$
\End_K(X) \otimes \Z/2\Z \hookrightarrow \End(X_2)
$$
is also an injection (see  \cite[p. 177]{MumAV}). The image of this
homomorphism lies in $\End_{\Gal(K)}(X_2)$.

Before we continue, let us prove the following useful statement,
which explains why the discriminants of orders considered in 
this paper have to belong to a certain congruence class.
\begin{lemma}\label{modlem}
Let $D$ be a $\Q$-algebra that is a $2$-dimensional
vector space over $\Q$, and let $\O$ be an order of $D$.
Then $\F_2$-algebras $\O\otimes\Z/2\Z$ and $\F_4$ are isomorphic
if and only if $D$ is a quadratic field and $\discr(\O)\equiv5\bmod8$.
\end{lemma}
\begin{proof}

First, we show that if $\O\otimes\Z/2\Z \cong \F_4$ then $D$ is a field.
Proceeding by contradiction, assume that for some nonzero $u, v \in D$, 
we have $uv = 0$.
After multiplying $u$ and $v$ by appropriate nonzero rational numbers
we can assume that $u, v \in \O$ and neither is divisible by $2$ in $\O$.
Then the image of neither $u$ nor $v$ in $\O\otimes\Z/2\Z$
equals to $0$, but their product does, which leads to a contradiction.
Hence $D$ is a $\Q$-division algebra of dimension $2$
as a $\Q$-vector space, or, equivalently, a quadratic field.

Assume $D$ is a quadratic field.
Let $d\geq 2$ be the reduced discriminant of $D$;
that is, assume that $D=\Q(\sqrt{d})$ with $d$ square-free.
It is well-known that $\O=\Z[c\eta]$, where
$$
\eta =
\begin{cases}
(\sqrt{d}-1)/2 & \text{for }d\equiv 1\bmod 4,\\
\sqrt{d} & \text{for }d\equiv 2\text{ or }3\bmod 4,\\
\end{cases}
$$
and $c$ is a positive integer called the {\em conductor} of $\O$.
The minimal polynomial of the generator $c\eta$ is
$$
g(X) =
\begin{cases}
X^2+cX-c^2(d-1)/4 & \text{for }d\equiv 1\bmod 4,\\
X^2-c^2d & \text{for }d\equiv 2\text{ or }3\bmod 4.\\
\end{cases}
$$
The field $\F_4$ is the splitting field
of $X^2+X+1$ over $\F_2$, 
the only irreducible polynomial of degree $2$ over $\F_2$.
We have $\Z[c\eta]\otimes\Z/2\Z\cong\F_4$
if and only if the minimal polynomial of
$c\eta$ is congruent to $X^2+X+1$ modulo $2$.
This happens if and only if $d\equiv5\bmod8$ and $c\equiv1\bmod2$.
It is easy to see that this is equivalent to
$\discr(\O)\equiv5\bmod8$.
\end{proof}

As an application of this lemma we obtain:
\begin{proof}[Proof of Theorem \ref{rankthm}]
Let $K$, $f$, $C_f$, and $X=J(C_f)$ be as in Theorem \ref{rankthm}.
By Theorem \ref{auxthm}, $\End_{\Gal(K)}(X_2) \cong \F_4$.
Hence $\End_K(X) \otimes \Z/2\Z \subset \F_4$, and
$\rank_{\Z} (\End_K(X)) \leq 2$. Rank of the
free $\Z$-module $\End_K(X)$ and dimension of $\Q$-algebra $\End_K^0(X)$
are equal to the $\F_2$-dimension of the algebra $\End_K(X)\otimes\Z/2\Z$.
So if $\rank_{\Z} (\End_K(X))=1$ then $\End_K(X)) = \Z$.

If $\rank_{\Z} (\End_K(X))=2$ then the $\Q$-algebra
$\End_K^0(X)$ has dimension $2$ as a $\Q$-vector space.
It is well-known that $\End_K(X)$ is isomorphic to an order of 
$\End_K^0(X)$ and by Lemma \ref{modlem}, we have 
$D$ is a quadratic field and
$$
\discr(\End_K(X)) \equiv 5 \bmod 8.
$$
\end{proof}

%----------------------------------------------
\begin{proof}[Proof of Theorem \ref{mainthm}]
Let $C=C_f$ be the hyperelliptic curve defined over $K$ by the equation $y^2=f(x)$
such that all conditions of Theorem \ref{mainthm} hold for the polynomial $f$ 
and abelian variety $X=J(C)$.

%According to \cite[p. 177]{MumAV}, if a subring $\O$ of $\End(X)$ satisfies
%$$
%(\O \otimes \Q) \cap \End(X) = \O
%$$
%then there exists an injective $\Gal(K)$-homomorphism
%$$
%\O \otimes \Z/2\Z \hookrightarrow \End(X_2).
%$$
%In particular, this is true if $\O$ equals to an intersection
%of $\End(X)$ and a subalgebra of $\End^0(X)$.

% ???? Fix this

We know that $\End_K(X)$ contains the order $\i(D) \cap \End(X)$
of $\i(D)$, whose $\Z$-rank is $2$. Therefore,
$\rank_{\Z} (\End_K(X)) = 2$ and, by Theorem \ref{rankthm},
$\End_K(X)$ is isomorphic to an order of $D$.

Let us consider the possible options for the $\Gal(K)$-stable algebra
$$
R:=\End(X)\otimes\Z/2\Z\subset\End_{\F_2}(X_2)
$$
which are provided by Theorem \ref{auxthm}. First, note that the rank of the
free $\Z$-module $\End(X)$ and dimension of $\Q$-algebra $\End^0(X)$
are equal to the $\F_2$-dimension of the algebra $\End(X)\otimes\Z/2\Z$.

\noindent{\bf Case (\ref{csf2}): $\End(X)\otimes\Z/2\Z \cong \F_2$.}
This case cannot occur, since the rank of $\End(X)$, which contains $\End_K(X)$,
is at least $2$.

\noindent{\bf Case (\ref{csf4}): $\End(X)\otimes\Z/2\Z \cong \F_4$.}
In this case, the free $\Z$-module $\End(X)$ has rank $2$,
and the ring $\End(X)$ is isomorphic to an order 
of the real quadratic field $D\cong\End_K^0(X)=\End_K(X)\otimes\Q$.
Note that $\End^0(X)\cong D$ is a simple division algebra,
so $X$ is not isogenous over the algebraic closure $K_a$ of $K$
to a product of two elliptic curves.
Therefore, $X$ is an absolutely simple abelian variety.
In this case, the equality $\End(J(C_f))=\End_K(J(C_f))$ holds.

\noindent{\bf Case (\ref{csm2f4}): $\End(X)\otimes\Z/2\Z\cong\Mat_2(\F_4)$.}
Then we have $\dim_{\Q}(\End^0(X))=8$.
In order to eliminate this outcome, let us consider the following possibilities.
\begin{enumerate}
\item
It is well-known \cite{OortEnd} that if $X$ is an absolutely simple abelian 
variety of dimension $2$, then its endomorphism algebra $\End^0(X)$ is
an Albert algebra of type I($1$), I($2$), II($1$), or IV($2$, $1$),
which means that $\dim_\Q(\End^0(X))\neq 8$.

\item
Suppose $X$ is isogenous over $K_a$ to a product of two non-isogenous elliptic curves
$E_1$ and $E_2$. We have
$$
\End^0(X)\cong\End^0(E_1) \oplus \End^0(E_2)
$$
and
$$
\dim_{\Q}(\End^0(X))=
\dim_{\Q}(\End^0(E_1))+
\dim_{\Q}(\End^0(E_2)).
$$
It is well-known \cite[pp.~102 and 165]{SilEllCurves},
that the endomorphism algebra of an elliptic curve has 
$\Q$-dimension of $1$ or $2$ in characteristic $0$, and $1, 2$ or $4$ in 
positive characteristic.
This means that $\dim_{\Q}(\End^0(X))\neq 8$,
since the equality would imply that $E_1$ and $E_2$ are 
both supersingular and, therefore, isogenous.

\item
If $X$ is isogenous over $K_a$ to a square of an
elliptic curve $E_1$ then
$$
\End^0(X)\cong\Mat_2(\End^0(E_1))
$$ 
and
$$
\dim_{\Q}(\End^0(X))=
4\dim_{\Q}(\End^0(E_1)).
$$
and for this dimension over $\Q$
to be $8$, we must have
$\dim_{\Q}(\End^0(E_1))=2$. This means that
$\End^0(X)$ is a matrix algebra of size $2$ over an imaginary
quadratic extension $L=\End^0(E_1)$ of $\Q$.
The order $\Cal L= L\cap\End(X)$ of the center $L$ 
of $\End^0(X)$ has the property
that $\Cal L\otimes\Z/2\Z\subset \End_{\F_2}(X_2)$ is
stable under the adjoint action of the group $\Gal(K)$.

We know that the ring $\End_K(X)$ has the same property.
Note that the subalgebra $\End_K^0(X) = \End_K(X) \otimes\Q$
of $\End^0(X)$ is isomorphic to a real quadratic field,
and therefore it does not coincide with the algebra $L$,
to which it is not isomorphic.
Hence the compositum $\End_K^0(X)L$ of $\End_K^0(X)$ and $L$ has
dimension $4$ over $\Q$. It is $\Gal(K)$-stable.

Let $\Cal R=\End_K^0(X)L \cap \End(X)$. 
The ring $\Cal R$ is an order of $\End_K^0(X)L$
and therefore has rank $4$ over $\Z$.
According to \cite[p. 177]{MumAV}, 
there exists an injective $\Gal(K)$-homomorphism
$\Cal R \otimes \Z/2\Z \hookrightarrow \End(X_2)$.
The image of this homomorphism is also $\Gal(K)$-stable
and therefore satisfies all of the conditions of 
Theorem \ref{auxthm}. 
Hence it must fit one of the three possible choices
prescribed by this theorem.
However, the $\F_2$-dimension of $\Cal R \otimes \Z/2\Z$ is $4$, 
while the dimensions of spaces in Theorem \ref{auxthm} 
are $1$, $2$, $8$, and $16$.
We arrive at a contradiction.
\end{enumerate}

\noindent{\bf Case (\ref{csm4f2}): $\End(X)\otimes\Z/2\Z=\End_{\F_2}(X_2)$.}
The free $\Z$-module $\End(X)$
has rank $16$. Recall that $g=\dim(X)=2$.
This implies $\rank_\Z  (\End(X))= (2g)^2$,
and the semisimple $\Q$-algebra
$\End^0(X)=\End(X)\otimes\Q$ has dimension $(2g)^2$.
This means that $\chr(K)>0$ and $X$ is a 
supersingular abelian variety (see \cite[Lemma~3.1]{ZarCM}).
\end{proof}

%===================================================================
\section{Permutation Groups and Permutation Modules}
%\label{auxsec2}
%===================================================================
In proving Theorem \ref{auxthm} we will need some basic results about
the structure of $J(C)_2$.
Let $B$ be a set of even cardinality $n\geq 6$.
Denote by $\Perm(B)$ the group of permutations on $B$. 
A choice of an ordering on $B$ induces an isomorphism $\Perm(B)\cong\bS_n$. 
Let $\F$ be a field of characteristic $2$, 
and denote by $\F^B$ the $n$-dimensional $\F$-vector space of maps from $B$ to $\F$.
The action of $\Perm(B)$ on $B$ extends to an action of $\Perm(B)$
on $\F^B$ as follows: an element $\sigma\in\Perm(B)$ sends
a map $f:B\to\F$ to map $\sigma f:b\mapsto f(\sigma^{-1}(b))$.
The subspace 
$$
(\F^B)^0=\sset{f: B\to \F \mid \sum_{b\in B} f(b) = 0}.
$$
of $\F^B$ is stable under the action of $\Perm(B)$. 
In turn, $\Perm(B)$-module $(\F^B)^0$ contains a stable submodule
$\F \cdot 1_B$ of constant functions $B\to\F$.
Given a subgroup $G$ of $\Perm(B)$, we define the {\em heart} of 
the permutation representation of $G$ on $B$ over $\F$ to be
the quotient
$$
(\F^B)^{00}=(\F^B)^0/(\F\cdot 1_B).
$$
It is easy to show that $(\F^B)^{00}$ is a faithful $G$-module.

When $\F=\F_2$ we will write $Q_B$ instead of $(\F_2^B)^{00}$.
In this case, $Q_B$ can also be described as the set of equivalence 
classes of subsets of $B$ of even cardinality with symmetric difference
as sum where subsets complementary in $B$ are identified.

\begin{lemma}
Let $G$ be a subgroup of $\Perm(B)$.
Then we have an $\F_4[G]$-module isomorphism
$$
(\F_4^B)^{00} \cong Q_B \otimes_{\F_2} \F_4
$$
and an $\F_4$-algebra isomorphism
$$
\End_{G, \F_4}((\F_4^B)^{00}) \cong \End_G(Q_B) \otimes_{\F_2} \F_4.
$$
\end{lemma}
\begin{proof}
The first statement is obvious and the second immediately
follows from Lemma 10.37 of \cite{CuMRT}.
\end{proof}

Let $C_f: y^2 = f(x)$ be a hyperelliptic curve defined over a field $K$ of characteristic
different from $2$ by an irreducible separable polynomial $f(x)\in K[x]$ of even degree $n$,
and let $\fR$ denote the set of roots of $f$.
It is well-known that the $\Gal(K)$-modules $J(C_f)_2$ and $Q_{\fR}$ are isomorphic.

Now assume that $n=6$ and $\Gal(f)\cong\PSL_2(\F_5)$.
Then $\dim(J(C_f))=2$ and $\dim_{\F_2}(J(C_f)_2)=4$.
Let $\orhox: \Gal(K) \to \Aut(J(C_f)_2)$ be the action of $\Gal(K)$ on $J(C_f)_2$ and
let $G = \orhox(\Gal(K)) \subset \Aut(J(C_f)_2)$ be the image of $\Gal(K)$ under this representation.
The action of $\Gal(K)$ on $J(C_f)_2$ factors through $G$, and the action
of $\Gal(K)$ on $\fR$ factors through $\Gal(f)$.
We have $G \cong \Gal(f) \cong \PSL_2(\F_5)$ and
the faithful $G$-modules $J(C_f)_2$ and $Q_\fR$ are isomorphic.

%--------------------------------------------------------------------------
\begin{lemma}\label{lem1}
$Q_\fR$ is a simple $G$-module.
\end{lemma}
\begin{proof}
See Table 1 in \cite{MortMPR}.
%(see also Theorem 3.3\cite{IvaFATG}).
\end{proof}
%--------------------------------------------------------------------------
\begin{lemma}\label{lem2} 
The $\F_2$-algebras $\End_G(Q_\fR)$ and $\F_4$ are isomorphic.
\end{lemma}
\begin{proof}
Since the representation $G \to \Aut_{\F_2}(Q_\fR)$ is irreducible,
by Schur's Lemma $\End_G(Q_\fR)$ is a division algebra. Since it is
finite, it must be a field, which we will denote by $\F$.

According to  Table 1 in \cite{MortMPR},
the $G$-module $(\F_4^B)^{00}$ is reducible. Therefore,
$\End_{G,\F_4}((\F_4^B)^{00}) \cong \End_G(Q_\fR) \otimes_{\F_2} \F_4$
is not a field. 
This means that $\F=\End_G(Q_\fR)$ contains $\F_4$,
a quadratic extension of $\F_2$, as a subfield.
Hence $\F\cong\F_{4^s}$ for some positive integer $s$.
The embedding of $\F$ in $\End_{\F_2}(Q_\fR)$ provides $Q_\fR$
with a structure of an $\F$-vector space. Since $\#(Q_\fR)=16$, we must have
$\F\cong \F_4$ or $\F_{16}$. 

If $\F \cong \F_{16}$ then $Q_\fR$ is a $1$-dimensional $\F$-vector
space, so $\End_{\F}(Q_\fR)=\F$. 
Since $\F=\End_G(Q_\fR)$, we have $G\subset\Aut_{\F}(Q_\fR)$ as a subgroup,
which is a contradiction. Indeed, the group $\Aut_{\F}(Q_\fR) = \F\mult$
is abelian, while $G\cong\A_5$ is not.
\end{proof}

In light of the fact that
$\End_G(Q_\fR)=\End_{\Gal(f)}(J(C)_2)=\End_{\Gal(K)}(J(C)_2)$
this lemma can be restated as
$$
\dim_{\F_2} \End_{\Gal(K)}(J(C)_2)=2,
$$
and, since $\End_K(J(C))\otimes\Z/2\Z$ injectively embeds in 
$\End_{\Gal(K)}(J(C)_2)$, we get Theorem \ref{rankthm}.

%--------------------------------------------------------------------------
From Lemmas \ref{lem1} and \ref{lem2}, and 
Theorem 3.43 of \cite[p. 54]{CuMRT} follows:
\begin{cor}\label{abssmplite}
The $\F_4[G]$-module $Q_\fR$ is absolutely simple.
\end{cor}

%===================================================================
\section{Proof of the Auxiliary Theorem}
\label{auxpf}
%===================================================================
Theorem \ref{auxthm} follows immediately from the discussion
of  the previous section and this theorem:
\begin{thm}\label{alglist}
Let $X$ be an abelian variety over field $K$ of characteristic different from $2$.
Let $G = \orhox(\Gal(K)) \subset \End_{\F_2}(X_2)$. Assume that the following
conditions are satisfied:
\renewcommand{\theenumi}{\alph{enumi}}
\begin{enumerate}
\item $\dim(X) = 2$,
\item $G\cong \A_5$,
\item $X_2$ is a simple $G$-module,
\item $\End_G(X_2) \cong \F_4$.
\end{enumerate}
Identify $\F_2$ with its embedding $\F_2\cdot\Id_{X_2}\subset \End(X_2)$, 
where $\Id_{X_2}$ is the identity automorphism of $X_2$,
and identify $\F_4$ with $\End_G(X_2)$.
Let $R$ be a subalgebra of $\End_{\F_2}(X_2)$ containing the
identity automorphism $\Id$ of $X_2$ such that
$$
uRu\inv\subset R\quad\text{for all }u\in G.
$$
Then we have one of the following cases:
\renewcommand{\theenumi}{\roman{enumi}}
\begin{enumerate}
\item $R = \F_2$;
\item $R = \F_4$;
\item $R = \End_{\F_4}(X_2)$;
\item $R = \End_{\F_2}(X_2)$.
\end{enumerate}
\end{thm}

\begin{proof}
Since $X_2$ is a faithful $R$-module, we have
$$
uRu\inv=R\quad\text{for all }u\in G\subset\Aut_{\F_2}(X_2).
$$
%--------------------------------------------------------------------------
\begin{lemma}
$X_2$ is a semisimple $R$-module.
\end{lemma}
\begin{proof}
{\em This is a reproduction of 
a similar proof in \cite[of Th.~5.3]{ZarCM}.}
Let $U\in X_2$ be a simple $R$-submodule.
Then $U'=\sum_{s\in G}sU$ is a
non-zero $G$-invariant subspace in $X_2$, 
and, since $X_2$ is a simple $G$-module, 
$U'=X_2$.
Each $sU$ is also an $R$-submodule in $X_2$, 
because $s^{-1}Rs=R$ for all $s\in G$. 
In addition, if $W\subset sU$ is an $R$-submodule
then $s^{-1}W$ is an $R$-submodule in $U$, because
\[
Rs^{-1}W=s^{-1}sRs^{-1}W=s^{-1}RW=s^{-1}W.
\]
Since $U$ is simple, $s^{-1}W=\set{0}$ or $U$. 
This implies that $sU$ is also simple. 
Hence $X_2=U'$ is a sum of simple
$R$-modules and therefore is a semisimple 
$R$-module.
\end{proof}
%--------------------------------------------------------------------------

\begin{lemma}
The $R$-module $X_2$ is isotypic.
\end{lemma}
\begin{proof}
{\em The proof is a modification of a similar proof \cite[of Th.~5.3]{ZarCM}.}
Let
$$
X_2 = V_1 \oplus \cdots \oplus V_r
$$
be an isotypic decomposition of the 
semisimple $R$-module $X_2$. 
Looking at the dimensions yields $r\leq \dim_{\F_2}(X_2)=4$. 
By repeating the argument in the proof of the previous 
claim, we can show that for each isotypic component $V_i$
its image $sV_i$ is an isotypic $R$-submodule for
each $s\in G$ and therefore is contained in some $V_j$.
Similarly, $s\inv V_j$ is an isotypic submodule containing $V_i$.
Since $V_i$ is the isotypic component, $s\inv V_j=V_i$.
This means that $s$ permutes the $V_i$, and, 
since $X_2$ is $G$-simple, $G$ permutes them transitively. 
This gives a homomorphism $G\to\bS_r$ which must 
be injective or trivial, since $G$ is simple.
However $G\cong\A_5$ and $r\leq 4$, so it is trivial. 
This means that $sV_i=V_i$ for all $s\in G$ and $X_2=V_i$ 
is isotypic. 
\end{proof}

From this lemma it follows that there exists a simple $R$-module $W$ and 
a positive integer $d$ such that $X_2\cong W^d$.

We have
$$
d \cdot \dim_{\F_2}(W) = \dim_{\F_2}(X_2) = 4.
$$
Thus $d=1, 2,$ or $4$. 

Clearly, $\End_R(X_2)$ is isomorphic to 
the matrix algebra $\Mat_d(\End_R(W))$.
Let us put
$$
k=\End_R(W).
$$
Since $W$ is simple, $k$ is a finite 
division algebra of characteristic $2$.
Hence $k$ is a finite field of characteristic $2$, and
$$
\End_R(X_2)\cong\Mat_d(k).
$$
We have $\End_R(X_2)\subset\End_{\F_2}(X_2)$ is invariant 
under the adjoint action of $G$, since $R$ is invariant
under adjoint action of $G$. 
This induces a homomorphism
$$
\alpha:G\to\Aut(\End_R(X_2))=\Aut(\Mat_d(k)).
$$
Since $k$ is the center of $\Mat_d(k)$, 
it is invariant under the action of $G$;
that is, we get a homomorphism $G\to\Aut(k)$, 
which must be trivial, since $G$ is a simple 
group and $\Aut(k)$ is abelian. 
This implies that the center $k$ of $\End_R(X_2)$ 
commutes with $G$ and must be
a subalgebra of $\End_G(X_2)$.
Since $\End_G(X_2) = \F_4$ as an $\F_2$-algebra, 
we have $k=\F_2$ or $\F_4$.

It follows from the Jacobson density theorem 
(combined with dimension arguments) that 
$R\cong\Mat_m(k)$ with $dm=4$ if $k=\F_2$
and $2dm=4$ if $k=\F_4$. 

Let us rule out the case not mentioned in the outcomes of this theorem:
if $R\cong\Mat_2(\F_2)$, the group $G$ acts on $\Mat_2(\F_2)$;
that is, we have a homomorphism
$$
G \to \Aut(\Mat_2(\F_2)) = \PGL_2(\F_2) \cong \GL_2(\F_2),
$$
where the equality follows from the Skolem-Noether theorem.
This homomorphism must be trivial, since $G$ is perfect and 
$\GL_2(\F_2)$ is solvable.
Therefore, $R$ commutes with $G$ and is a subalgebra
of $\End_G(X_2)=\F_4$.
This is a contradiction.

Finally, notice that if $k=\End_G(X_2) = \F_4$, and
$R=\Mat_2(k)$, then $R$ commutes with $\F_4$, and therefore
lies in $\End_{\F_4}(X_2)$. However, since $X_2$ is a $2$-dimensional
$\F_4$-vector space, $\End_{\F_4}(X_2)\cong \Mat_2(\F_4)$, so
$R=\End_{\F_4}(X_2)$.
\end{proof}

%-----------------------------
In \cite{ZarMR}, Yu.~Zarhin makes the following definition.
\begin{defn}
Let $V$ be a vector space over a field $k$, let $G$ be a group and 
$\rho: G \to \Aut_k(V)$ a linear representation of $G$ in $V$. We say that
the $G$-module $V$ is {\em very simple} if it enjoys the following property:

If $R \subset \End_k(V)$ is a $k$-subalgebra containing the identity operator $\Id$
such that
$$
\rho(\sigma)R\rho(\sigma)\inv \subset R \quad \text{for all }\sigma\in G
$$
then either $R=k\cdot\Id_V$ or $R=\End_k(V)$.
\end{defn}

%-----------------------------
We immediately obtain the following proposition, which will be used in
the proof of Theorem \ref{nsauxthm}.
\begin{cor}\label{verysimplecor}
Let $X$ be an abelian variety satisfying the conditions of Theorem \ref{alglist}.
Then $X_2$ is a very simple $G$-module over $\F_4$.
\end{cor}
\begin{proof}
The only outcomes of Theorem \ref{alglist} that contain $\F_4$
in the center are $R \cong \F_4$ and $R =\End_{\F_4}(X_2)$.
\end{proof}

%===================================================================
\section{Non-Isogenous Jacobians}
\label{nonisogenuity}
%===================================================================
Let us now restate some of the above results in terms of $J(C)_2$ and $\Gal(f)$.
If we compose the canonical epimorphism $\Gal(K)\twoheadrightarrow\Gal(f)$ with 
the irreducible representation $\Gal(f)\to\Aut_{\F_2}(J(C)_2)$ of Lemma \ref{lem1} 
we get
\begin{lemma}\label{lem1b}
$J(C)_2$ is a simple $\Gal(K)$-module.
\end{lemma}

Note that Corollary \ref{abssmplite} can be restated as follows:
\begin{cor}\label{abssmp}
Assume that all of the conditions of Theorem \ref{mainthm} hold.
Then $J(C_f)_2$ is an absolutely simple $\F_4[\Gal(f)]$-module.
\end{cor}

\begin{proof}[Proof of Theorem \ref{NICor}]
We prove this theorem by contradiction. Let $X=J(C_f)$ and $Y=J(C_h)$ be abelian
varieties in question and assume there exists a nonzero homomorphism $\phi\in\Hom(X, Y)$.
Then $\phi$ is an isogeny. 
%Indeed, the connected component of the kernel of $\phi$ containing
%zero is an abelian subvariety of $X$ and so is the image of $X$ under $\phi$. 
%Since both $X$ and $Y$ are simple and $\phi$ is not identically zero, 
%the kernel is finite and the image is the whole of $Y$.

We can also assume that $X_2\not\subseteq\ker\phi$. If that is not the case,
then $\phi$ is a composition of multiplication by $2$ on $X$ and another isogeny 
from $X$ to $Y$, which we can choose instead, continuing this process until
the obtained isogeny no longer annihilates $X_2$.

For every $\psi\in\Hom(X, Y)$ and $\sigma\in\Gal(K)$ we define a homomorphism
${^\sigma}\!\psi \in \Hom(X, Y)$ by ${^\sigma}\!\psi (x) = \sigma\psi(\sigma\inv x)$ for all $x\in X(K_a)$.
We then define $c:\Gal(K)\to\End^0(X)\mult$ by 
$c_\sigma=\phi\inv\,{^\sigma}\!\phi$. 
It is easy to show that $c$ satisfies the cocycle condition 
$c_{\sigma\tau}=c_\sigma\,{^\sigma}\!c_\tau$. 
In addition,  we have $\End^0(X)=\End_K^0(X)$,
so ${^\sigma}\!c_\tau=c_\tau$ and $c$ is a homomorphism.
There exists a finite normal extension $K_\phi$ of $K$ such that 
$\phi$ is defined over $K_\phi$. 
Since we are in characteristic $0$, 
$K_\phi$ is also Galois over $K$, so there exists a homomorphism 
$c':\Gal(K_\phi/K)\to\End^0(X)\mult$ such that $c$ is a composition of the
canonical homomorphism $\Gal(K)\twoheadrightarrow\Gal(K_\phi/K)$
and $c'$. Since $\Gal(K_\phi/K)$ is finite, its image under $c'$ in $\End^0(X)\mult$,
which coincides with the image of $\Gal(K)$ under $c$, is also finite.
Therefore, this image is either $\set{\Id_X}$ or $\mu_2=\set{\pm\Id_X}$.

If $c(\Gal(K))=c'(\Gal(K_\phi/K))=\mu_2$, let $H=\ker(c')$,
and let $K'$ be the subfield of elements of $K_\phi$ that are fixed by $H$. 
Then $K'$ is Galois over $K$,
$\Gal(K_\phi/K')=H$ and $\Gal(K'/K)\cong G/H=\mu_2$.
By the choice of $H$, the image of $\Gal(K_\phi/K')$ under $c'$ is trivial.
Therefore, the image of $\Gal(K') \subset \Gal(K)$ under $c$ is trivial.
Moreover, since the Galois groups of $f$ and $h$ are perfect and $\Gal(K'/K)\cong\mu_2$
is cyclic, polynomials $f$ and $h$ will be irreducible over $K'$, their
Galois groups over $K'$ will still be isomorphic to $\A_5$, and
they will still satisfy the conditions of Theorem \ref{NICor}.
Without loss of generality, we can choose to work over $K'$ instead of $K$, 
which reduces the theorem to the next case.

If $c(\Gal(K))=\set{\Id_X}$ then $\phi$ is defined over $K$ and 
commutes with the action of $\Gal(K)$ on $X$ and $Y$. 
This remains true if we consider points of order dividing $2$:
the homomorphism $\varphi = \phi|_{X_2}: X_2 \to Y_2$ commutes with the action
of $\Gal(K)$ on $X_2$ and $Y_2$.
The kernel of $\varphi$ is a $\Gal(K)$-stable
submodule of $X_2$. Since $X_2$ is a simple $\Gal(K)$-module, 
the map $\varphi$ is either zero or is a $\Gal(K)$-isomorphism.
It cannot be zero by the choice of $\phi$.
We claim that it cannot be an isomorphism either.

Let $L=K(\fR_f\cup \fR_h)$ be the compositum of the splitting fields
$K(\fR_f)$ and $K(\fR_h)$ of $f$ and $h$ and consider 
the canonical restriction maps $\res_{K(\fR_f)}^L:\Gal(L/K) \to \Gal(f/K)$
and $\res_{K(\fR_h)}^L:\Gal(L/K) \to \Gal(h/K)$. 
Recall that
the $\Gal(f/K)$-module $X_2$ and $\Gal(h/K)$-module $Y_2$
are faithful and
let $\alpha_X:\Gal(f/K) \hookrightarrow \Aut(X_2)$
and $\alpha_Y:\Gal(h/K) \hookrightarrow \Aut(X_2)$ be
the corresponding embeddings.
The actions of $\Gal(K)$ on $X_2$ and $Y_2$ factor through
the canonical epimorphism $\res_L: \Gal(K)\twoheadrightarrow\Gal(L/K)$,
that is,
$$
\orhox = \alpha_X \circ \res_{K(\fR_f)}^L \circ \res_L
\quad\text{ and }\quad
\orhoy = \alpha_Y \circ \res_{K(\fR_h)}^L \circ \res_L.
$$

If $K(\fR_f)$ and $K(\fR_h)$ are linearly disjoint over $K$, then
$$
\Gal(L/K) \cong \Gal(f/K) \times \Gal(h/K),
$$
with projections onto each summand coinciding with the Galois
restriction maps $\res_{K(\fR_f)}^L$ and $\res_{K(\fR_h)}^L$.
Let $\sigma_L = (\sigma_f, \Id_h)\in \Gal(f/K) \times \Gal(h/K) = \Gal(L/K)$, 
where $\Id_h$ is the identity element of $\Gal(h/K)$ and
$\sigma_f \in \Gal(f/K)$ is a nonidentity element,
and pick any $\sigma \in  \res_L\inv(\sigma_L)$.
Then
$$
\orhox(\sigma) = \alpha_X(\res_{K(\fR_f)}^L(\sigma_L)) = \alpha_X(\sigma_f)\neq \Id_{X_2},
$$
while
$$
\orhoy(\sigma) = \alpha_Y(\res_{K(\fR_h)}^L(\sigma_L)) = \alpha_Y(\Id_h) = \Id_{Y_2}.
$$
Therefore, the $\Gal(K)$-modules $X_2$ and $Y_2$ are not isomorphic.

\end{proof}

%====================================================================
\section{The centralizer of endomorphisms defined over base field.}
\label{centr}
%====================================================================
This section contains preliminary investigations required for the proof
of Theorem \ref{nsthm} about the impossibility of the supersingular outcome in
characteristics $p>2$ when $p$ splits in the quadratic field $D = \Q(\omega) \cong \End_K^0(X)$. 
First, we determine the algebraic structure and places of ramification of the centralizer of 
$\i(D)=\End_K^0(X)$ in $\End^0(X)$ in the case when $X$ is supersingular.
Put
\begin{eqnarray*}
\End^0(X, \i) &=&
\sset{u\in \End^0(X)
\mid
\i(y)u = u\i(y) \quad \forall y\in D}\\
&=&
\sset{u\in \End^0(X)
\mid
\i(\omega)u = u\i(\omega)}.
\end{eqnarray*}
It is well-known that when $X$ is a supersingular abelian surface,
then $\End^0(X)\cong\Mat_2(\H_p)$,
where $\H_p$ a quaternion $\Q$-algebra ramified exactly at $p$ and $\infty$. 
We write $\H_{p, D}$ for $\H_p \qtimes D$. Both $\End^0(X, i)$ and $\H_{p, D}$
carry natural structures of $D$-algebras.

\begin{thm}
The $D$-algebras $\End^0(X, \i)$ and $\H_{p, D}$ are isomorphic.
\end{thm}

\begin{proof}
We begin by showing that the isomorphism class of the $D$-algebra $\End^0(X, \i)$ 
is independent of the embedding of $D$ into $\End^0(X)$
that sends $1$ to the identity automorphism of $X$.
Let $j: D\hookrightarrow\End^0(X)$ be another such embedding.
The $\Q$-algebra $\i(D)$ is a simple $\Q$-subalgebra of the simple central $\Q$-algebra
$\End^0(X)$ and $j\i\inv:\i(D)\to j(D)$ is an algebra isomorphism.
By Skolem-Noether theorem \cite[p. 69]{CuMRT}, 
there exists $\sigma\in \Aut(\End^0(X))$
such that $ji\inv(x)=\sigma(x)$ for all $x\in\i(D)$.
Then $\sigma$ is a $\Q$-algebra isomorphism between $\End^0(X, \i)$ and $\End^0(X, j)$.
Indeed, for every $z \in \End^0(X, \i)$, we have
$$
\i(\omega)\inv z \i(\omega) = z
$$
and so
$$
\sigma(\i(\omega))\inv \sigma(z) \sigma(\i(\omega)) = \sigma(z).
$$
Taking into account that $\sigma(\i(\omega))=j(\omega)$, we obtain
$$
j(\omega)\inv \sigma(z) j(\omega) = \sigma(z),
$$
that is, $\sigma(z) \in \End^0(X, j)$. The map $\sigma: \End^0(X, \i) \to \End^0(X, j)$ is
obviously invertible.

Fix an isomorphism $\End^0(X) \cong \Mat_2(\H_p)$ and
let $j: D \hookrightarrow\End^0(X)$ be given by
$$
1\mapsto\Id=\begin{pmatrix} 1 & 0 \\ 0 & 1\end{pmatrix}, \quad 
\omega\mapsto\begin{pmatrix} 0 & d \\ 1 & 0\end{pmatrix}.
$$
Assume $M = \begin{pmatrix} r & s \\ t & u\end{pmatrix}\in \Mat_2(\H_p)$ 
commutes with $j(\omega)$.
Then $r=u$ and $s=td$,
$$
M=
\begin{pmatrix}
r & td\\
t & r
\end{pmatrix}
=r \Id + t j(\omega),
$$
and the centralizer of $j(D)$ in $\End^0(X)$ is generated over $\H_p$
by $\Id$ and $j(\omega)$. Therefore $\End^0(X,j)$ (\ $\cong\End^0(X,\i)$\ ) is isomorphic to 
$\H_p + j(\omega) \H_p = j(D)\H_p$, the subalgebra generated by the products of 
elements in $j(D)$ and $\H_p$.

Let $\H_p \qtimes D \to D\H_p$ be the $\Q$-algebra homomorphism induced by
$a \otimes b \mapsto ab$.
Since $\H_p \qtimes D$ is a simple algebra, the kernel of this map
is either trivial or is the entire $\H_p \qtimes D$. Because this homomorphism is not zero,
it must be an embedding. A comparison of dimensions over $\Q$ yields
$\End^0(X, i) \cong D\H_p \cong \H_p \qtimes D$.
\end{proof}

It is well-known \cite[p. 4]{VigQA} that $\H_p \qtimes D$ is a quaternion algebra over $D$.
\begin{lemma}\label{quatlem}
Let $\p$ be a place of $D$ dividing $p$.
The quaternion algebra $\End^0(X, i)$ ramifies at
$\p$ if and only if $p$ splits in $D$,
\end{lemma}
\begin{proof}
For a prime $\p \divides p$ of $D$, 
the degree $[D_\p:\Q_p]$ equals to the product $e_\p f_\p$ of ramification index
and relative degree of $\p$ over $p$.
This product equals $1$ if $p$ splits in $D$ and $2$ otherwise.
{\em Cancellation of ramification} at $p$ occurs if and only if the degree of this extension is even
(see \cite[Ch.~II, Th.~1.3, p.~33]{VigQA}).
\end{proof}

%====================================================================
%\section{Schur indices over localizations}\label{schur}
%====================================================================
In Section \ref{NSSect} we will demonstrate that if $X$ is a supersingular abelian variety,
then $\End^0(X, \i)$ is isomorphic to a direct summand of the group algebra $D[\SL_2(\F_5)]$. 
The following lemma will be used in the proof of Theorem \ref{nsauxthm}:
\begin{lemma}\label{schurlem}
Let\; $\Gamma = \SL_2(\F_5)$, let $D$ be a real quadratic field,
and $p$ a prime such that $\H_{p,D}=\H_p\otimes_{\Q}D$ is a quaternion
$D$-algebra.
Suppose that the group algebra $D[\Gamma]$ contains a direct summand
isomorphic to $\H_{p, D}$. Then the rational prime $p$ does not split in $D$.
\end{lemma}

\begin{proof}
Since
$$
\End^0(X, \i) \dtimes \C \cong \Mat_2(\C),
$$
the simple algebra $\H_{p, D}$ corresponds to a faithful irreducible 
character $\chi$ of $\Gamma$ over $D$ of degree $2$.

In the notation of \cite[Table~II]{FeitSI}, $\chi=\zeta_i$, $i=1$ or $2$.
Over $\Q$ we have $m_\ell(\chi)=m_{\Q_\ell}(\chi)=1$ for all
places $\ell$ of $\Q$ except for $\infty$. % , over which the Schur index is $2$.
This result can be extended to primes of $D$ by means of Theorem 2.16
in \cite{FeitSI}, which says that for extension $D_\q$ of $\Ql$, where $\q$ is
a prime of $D$ diviging $\ell$, we have
$$
m_\q(\chi)=\frac{m_\ell(\chi)}{(m_\ell(\chi), [D_\q(\chi): \Ql(\chi)])}.
$$
Therefore,
$$
m_\q(\chi)=m_{D_\q}(\chi)=1.
$$
for $\q \divides \ell$. 
This means that $\H_{p, D}$ does not ramify at any non-archimedean place of $D$,
and, by Lemma \ref{quatlem}, $p$ does not split in $D$.
\end{proof}

%====================================================================
\section{Non-Supersingularity}
\label{NSSect}
%====================================================================
In this section we prove the following theorem.
\begin{thm}\label{nsthm}
Assume the conditions of Theorem \ref{mainthm} hold for 
a quadratic field $D$ of reduced discriminant $d$, polynomial $f$, 
and abelian variety $X=J(C_f)$.
In addition, assume that $p=\chr(K)>2$ splits in $D$.
Then $X$ is not a supersingular abelian variety.
\end{thm}

Recall that if $\ell$ is an odd prime, then $\ell \divides d$ means that $\ell$ ramifies in $D$,
$(d/\ell)=-1$ means that $\ell$ is inert in $D$, and $(d/\ell)=1$ means that $\ell$ 
splits in $D$. Since $d\equiv 5\pmod 8$, the rational prime $2$ is inert in $D$,
and $D_2:=D \qtimes \Q_2$ is a field.

The proof of Theorem \ref{nsthm} will require the following proposition,
whose proof will be given later in this section.
\begin{thm}\label{nsauxthm}
Let $F$ be a field of characteristic $p>2$ containing all $2$-power roots of unity.
Let $G=\PSL_2(\F_5)$ and $D$ be a real quadratic field with $\discr(D)\equiv 5 \bmod 8$.
Suppose $X$ is an abelian variety defined over $F$, 
and assume that the following conditions hold:
\renewcommand{\theenumi}{\alph{enumi}}
\begin{enumerate}
\item $X$ is supersingular and $\dim(X)=2$;
\item there exists an injective $\Q$-algebra homomorphism 
$\i: D \hookrightarrow \End_F^0(X)$ such that
$\i(1)=\Id_X$, the identity automorphism of $X$;
\item the image of $\Gal(F)$ in $\Aut(X_2)$ is isomorphic to $G$ and the
corresponding faithful representation
$$
\orho: G \hookrightarrow\Aut(X_2)\cong\GL_4(\F_2)
$$
satisfies
$$
\End_G(X_2) = \F_4.
$$
\end{enumerate}
Then there exists a surjective group homomoprhism
$$
\pi_1:G_1\twoheadrightarrow G
$$
enjoying the following properties:
\renewcommand{\theenumi}{\roman{enumi}}
\begin{enumerate}
\item $G_1 \cong \SL_2(\F_5)$.
\item One may lift $\orho\pi_1: G_1\to\Aut(X_2)$ to a faithful absolutely irreducible
over the field $D_2 = D\qtimes\Q_2$ symplectic representation
$$
\rho: G_1 \hookrightarrow \Aut_{D_2}(V_2(X))
$$
in such a way that the following conditions hold:
\begin{itemize}
\item $\rho(G_1)\subset \End^0(X, \i)\mult$, where  $\End^0(X, \i)$ is the centralizer of
$\i(D)$ in $\End^0(X)$.
\item The homomorphism from the group algebra $D[G_1]$ to
$\End^0(X, \i)$ induced by $\rho$ is surjective and identifies
$\End^0(X, \i)$ with a direct summand of $D[G_1]$. % attached to $\chi$.
\end{itemize}
\end{enumerate}
\end{thm}

\begin{proof}[Proof of Theorem \ref{nsthm}]
Assume that $X=J(C_f)$ satisfies conditions of Theorem \ref{mainthm} and
is a supersingular abelian variety. 

Let $F\subset K_a$ be a field obtained from $K$ 
by adjoining all $2$-power roots of unity.
Then
$$
D \hookrightarrow \End_K^0(X) \subset \End_F^0(X).
$$
Moverover, the polynomial $f$ remains irreducible over $F$
and the Galois group of its splitting field over $F$ is still $\A_5$, a perfect group.
Indeed, $F$ is an abelian extension of $K$ and $\Gal(F/K)$ does not contain $\A_5$.
Thus $f$ and $X$ satisfy the conditions of Theorem \ref{mainthm} over $F$.
From Theorem \ref{nsauxthm} and Lemma \ref{schurlem}
it follows that $p$ does not split in $D$.
\end{proof}

\begin{proof}[Proof of Theorem \ref{nsauxthm}]
This proof is a modification of the proof of Theorem 3.3 in \cite{ZarNS}, and most of
what follows is actually stated in that work.

Let $T_2(X)$ denote the $\Z_2$-adic Tate module of $X$, 
$V_2(X):=T_2(X) \otimes_{\Z_2} \Q_2$ the $\Q_2$-adic Tate module of $X$,
and let
$$
\rhox: \Gal(F) \to \Aut_{\Z_2}(T_2(X))
$$
be the corresponding $2$-adic representation.
Put
$$
H:=\rhox(\Gal(F)) \subset \Aut_{\Z_2}(T_2(X)).
$$

%%%%%%%%%%%%%%%%%%%%%%%%%%%%%%%%%%%%%%%%%%%%%%%%
\begin{claim}\label{cl1}
The group $H$ is finite.
\end{claim}
\begin{proof}[Proof of Claim \ref{cl1}]
The rank of the $\Z_2$-module $T_2(X)$ is $2\dim(X)=4$, and, as Galois module,
$$
X_2 = T_2(X) / 2T_2(X).
$$
If we compose $\rhox$ with the surjective reduction modulo $2$ map
$$
\Aut_{\Z_2}(T_2(X)) \to \Aut(X_2),
$$
we get a natural homomorphism
$$
\orhox:\Gal(F) \to \Aut(X_2),
$$
which defines the action of $\Gal(F)$ on points of $X_2$.
Restriction of $\orhox$ to $H$ yields a natural continuous surjection
$$
\pi: H \to \orhox(\Gal(F))\cong G \subset \Aut(X_2).
$$

The choice of polarization on $X$ gives a non-degenerate alternating bilinear form (Riemann form)
$$
e:V_2(X) \times V_2(X) \to \Q_2.
$$
Since $F$ contains all $2$-power roots of unity, $e$ is $\Gal(F)$-invariant and hence $H$-invariant. 
Thus
$$
H\subset \Sp(V_2(X), e)
$$
and the $H$-module $V_2(X)$ is symplectic.

There exists a finite Galois extension $L$ of $F$ over which
all endomorphisms of $X$ are defined, that is,
$$
\End_L(X)=\End(X).
$$
Then then  group $\Gal(L)$ is an open subgroup of finite index in $\Gal(F)$, and the group
$$
H' := \rhox(\Gal(L))
$$
is an open normal subgroup of finite index in $H$.

There exists a natural embedding
$$
\End^0(X) \qtimes \Q_2 \hookrightarrow \End_{\Q_2}(V_2(X)),
$$
and, since $X$ is supersingular,
$$
\dim_\Q (\End^0(X)) = (2\dim(X))^2 = \dim_{\Q_2} (\End_{\Q_2}(V_2(X))),
$$
which implies that
$$
\End^0(X) \qtimes \Q_2 = \End_{\Q_2}(V_2(X)).
$$
Since all endomorphisms of $X$ are defined over $L$, the image $\rhox(\Gal(L))$ in
$\Aut_{\Q_2}(V_2(X))$ commutes with $\End^0(X)$, and therefore with the whole
$$
\End_{\Q_2}(V_2(X)) = \End^0(X) \qtimes \Q_2.
$$ 
This implies that 
$$
H'=\rhox(\Gal(L))\subset \Q_2\cdot \Id_{V_2(X)}.
$$
Since
$$
H'=\rhox(\Gal(L))\subset\rhox(\Gal(F)) \subset \Sp(V_2(X), e) \subset \SL(V_2(X)),
$$
we have
$$
H' \subset  \SL(V_2(X)) \cap \Q_2\cdot \Id_{V_2(X)},
$$
and the group $H' = \rhox(\Gal(L))$ is finite. 
As it is a subgroup of finite index in $H= \rhox(\Gal(F))$, 
the group $H$ is also finite.
\end{proof}

Since $H$ is finite, there exists a minimal subgroup $G_1$ of $H$ such that $\pi(G_1)=G$.
Denote the restriction of $\pi: H\to G$ to $G_1$ by $\pi_1:G_1 \twoheadrightarrow G$.
Put
\begin{equation}\label{defne}
E := \End_{G_1}(V_2(X)) \subset  \End_{\Q_2}(V_2(X)).
\end{equation}
%%%%%%%%%%%%%%%%%%%%%%%%%%%%%%%%%%%%%%%%%%%%%%%%
\begin{claim}\label{cl2}
The algebra $E$ is a quadratic field extension of $\Q_2$ and
$$
E \cong D \qtimes \Q_2.
$$
\end{claim}
\begin{proof}[Proof of Claim \ref{cl2}]
The $\Z_2$-algebra
$$
\O = E \cap  \End_{\Z_2}(T_2(X))
$$
is a free $\Z_2$-module, whose $\Z_2$-rank coincides with $\dim_{\Q_2}(E)$. 
The map
$$
\O/2\O \to \End_{\Z_2}(T_2(X)) / 2\End_{\Z_2}(T_2(X)) = \End(X_2)
$$
is an embedding. % justify
The $\Z_2$-rank of $\O$ equals the $\F_2$-dimension of the image of
$\O/2\O$ in $\End_G(X_2)$. % justify
Since elements of $\O$ commute with $G_1$ in $\End_{\Z_2}(T_2(X))$,
the image of $\O/2\O$ lies in $\End_G(X_2)\cong \F_4$. %justify
This implies
$$
\dim_{\Q_2}(E) = \rank_{\Z_2}\O = \dim_{\F_2}(\O/2\O) \leq 2.
$$

Let $F_1\subset F_a$ be the subfield fixed elementwise by
$$
\{\sigma \in \Gal(F) \mid \rhox(\sigma) \in G_1\}
$$
Then $F_1$ is a finite separable extension of $F$ and
$$
G_1 = \rhox(\Gal(F_1)).
$$
The image $\orhox(\Gal(F_1))$ in $\Aut(X_2)$ coincides with $G$.

Since $F\subset F_1$, 
we have $i(D) \subset \End_F^0(X) \subset \End_{F_1}^0(X)$ and so
\begin{eqnarray*}
\i(D) \qtimes \Q_2
&\subset& \End_{F_1}^0(X) \qtimes \Q_2\\
&\subset& \End_{\Gal(F_1)}(V_2(X))\\
& = & \End_{G_1}(V_2(X))\\
& = & E.
\end{eqnarray*}
Therefore
$$
\dim_{\Q_2}(\i(D) \qtimes \Q_2) = 2 \leq \dim_{\Q_2}(E)\leq 2,
$$
so $\dim_{\Q_2}(E)= 2$ and $E=\i(D)\qtimes \Q_2$.
We also get 
$$
\End_F^0(X)=\End_{F_1}^0(X)=\i(D).
$$
\end{proof}
%%%%%%%%%%%%%%%%%%%%%%%%%%%%%%%%%%%%%%%%%%%%%%%%

This provides $V_2(X)$ with a structure of a $2$-dimensional $E$-vector space
and gives us a faithful representation
$$
\rho:G_1\to\Aut_E(V_2(X)) \cong \GL_2(E),
$$
which must be absolutely irreducible by choice \eqref{defne} of $E$.

\begin{claim}\label{verysimpleclaim}
The $G_1$-module $V_2(X)$ is very simple over $E$.
\end{claim}
\begin{proof}[Proof of Claim \ref{verysimpleclaim}]

Since $\pi(G_1)=G$,
by \cite[Remark 5.2(i)]{ZarVS2R} and Corollary \ref{verysimplecor},
the $G_1$-module $X_2$ is very simple over $\F_4$.

Let $\O_E$ be the valuation ring of the quadratic $2$-adic field $E$,
and let $\lambda$ be its maximal ideal.
Since the prime $2$ of $\Q$ is inert in $D$, the prime ideal $\ell=2\Z_2$ of $\Z_2$
is also inert in $\O_E$, the completion of the ring of integers of $D$ with respect to
the $2$-adic topology. Hence $\lambda = \ell\O_E$.
Since the degree of inertia of $\lambda$ over $\ell$ is $2$, the residue field
$k(\lambda)=\O_E/\lambda$ is a quadratic extension of $\Z_2/\ell \cong \F_2$,
that is, $k(\lambda) \cong \F_4$.

The abelian group $T_2(X)$ is an $\O_E$-lattice in $V_2(X)$,
and therefore $T_2(X)/\lambda T_2(X)$ is a $k(\lambda)$-module.
Since $\lambda = \ell\O_E$, we have
$$
T_2(X)/\lambda T_2(X) = T_2(X)/(\lambda\O_E) T_2(X) = 
T_2(X)/ (\ell\O_E) T_2(X) = X_2.
$$
The $\O_E$-lattice $T_2(X)$ of $V_2(X)$ is $G_1$-stable,
so the $E[G_1]$-module $V_2(X)$ is a lifting of 
the very simple $\F_4[G_1]$-module $X_2$.
The claim follows from \cite[Remark 5.2(v)]{ZarVS2R}.
\end{proof}

%%%%%%%%%%%%%%%%%%%%%%%%%%%%%%%%%%%%%%%%%%%%%%%%
\begin{claim}\label{clc}
The group $G_1$ is a perfect central extension of $G$.
\end{claim}
\begin{proof}[Proof of Claim \ref{clc}]
Since $G \cong \PSL_2(\F_5)$ is perfect, so is $G_1$
(otherwise, we can replace $G_1$ with $[G_1, G_1]$, thus contradicting minimality of $G_1$).

By \cite[Remark 5.2(iv)]{ZarVS2R},
since the $G_1$-module $V_2(X)$ is very simple,
then either the normal subgroup $Z_1 = \ker (\pi_1: G_1 \to G)$ of 
$G_1$ consists of scalars (that is, it lies in $E$),
or the $E[Z_1]$-module $V_2(X)$ is absolutely simple.

We exclude the latter possibility by contradiction.
The kernel $Z$ of $\pi: H \to G$ is a subgroup of 
$1+2 \End_{\Z_2}(T_2(X)) \cong 1 + 2 \Mat_4(\Z_2)$.
Since $H$ is a finite group, so is $Z$.
In addition, $(1 + 2 \Mat_4(\Z_2))^2 \equiv 1 \pmod{ 2^2 \Mat_4(\Z_2)$}.
Thus, by Minkowski-Serre Lemma \cite[Th. 6.3]{SilMS},
the group $Z$ has exponent $1$ or $2$.

Therefore, $Z$ is a finite commutative group, and so is $Z_1 \subset Z$.
Hence $Z_1$ does not admit an absolutely irreducible
representation of dimension greater than $1$, which contradicts
$\dim_E(V_2(X))=2$.
Thus $Z_1 = \ker \pi_1 \subset E=\End_{G_1}(V_2(X))$ commutes with $G_1$.
\end{proof}

%%%%%%%%%%%%%%%%%%%%%%%%%%%%%%%%%%%%%%%%%%%%%%%%
\begin{claim}\label{cl3}
$G_1 \cong \SL_2(\F_5)$.
\end{claim}
\begin{proof}[Proof of Claim \ref{cl3}]
It is known \cite[Prop.~4.227 and Prop.~4.232(ii)]{GorFSG} that the only perfect central extensions of 
$\PSL_2(\F_5)$ are $\SL_2(\F_5)$ and $\PSL_2(\F_5)$ itself.
However, there are no faithful irreducible representations of $\PSL_2(\F_5)$ 
of dimension $2$ in characteristic $0$ (see \cite[p. 365]{CuMRT}).
\end{proof}

%%%%%%%%%%%%%%%%%%%%%%%%%%%%%%%%%%%%%%%%%%%%%%%%
\begin{claim}\label{cl4}
We have $\rho(G_1)\subset \End^0(X, i)\mult$.
\end{claim}
\begin{proof}[Proof of Claim \ref{cl4}]
We have $\orhox(\Gal(F_1))=G$ and hence 
$$
\End_{F_1}(X) \otimes \Z/2\Z \hookrightarrow \End_{\Gal(F_1)}(X_2)=\F_4
$$

Let $L_1$ be the finite Galois extension of $F_1$ attached to
$$
\rhox: \Gal(F_1) \to \Aut_{\Z_2}(T_2(X)).
$$
Then $\Gal(L_1/F_1)=G_1$. 
Since the image $\rhox(\Gal(L_1))$ in $\Aut_{\Z_2}(T_2(X))$ is trivial
and all $2$-power torsion points of $X$ are defined over $L_1$,
all endomorphism of $X$ are defined over $L_1$.
Hence there is a natural homomorphism
$$
\k : G_1 = \Gal(L_1 / F_1) \to \Aut(\End(X))
$$
such that
$$
\End_{F_1}(X)=\sset{u\in\End(X) 
\mid
\k(\sigma)u=u\quad\forall\sigma\in\Gal(L_1/F_1)}
$$
and
$$
\sigma(ux)=(\k(\sigma)u)(\sigma(x)).
$$
Further, we write $^{\k(\sigma)}\!u$ instead of $\k(\sigma)(u)$.
Since all $2$-power torsion points of $X$ are defined over $L_1$,
$$
\sigma(ux)= ^{\k(\sigma)}\!u(\sigma(x))
\text{ for all } x \in  T_2(X), u \in \End(X), \sigma \in G_1.
$$
Since $\Aut(\End(X)) \subset \Aut(\End^0(X))$, 
we can extend $\k$ to $\End^0(X)$.
Then
$$
\End_{F_1}^0(X)
=\{u\in\End^0(X) 
\mid
^{\k(\sigma)}\!u=u\quad\forall\sigma\in\Gal(L_1/F_1)\}
$$
and
$$
\sigma(ux)=^{\k(\sigma)}\!u(\sigma(x))
$$
Recall that
$$
\End^0(X) \subset \End^0(X) \qtimes \Q_2 = \End_{\Q_2}(V_2(X))
$$
and
$$
G_1 \subset \GL(V_2(X)) = \End_{\Q_2}(V_2(X))\mult.
$$
It follows that
$$
\sigma u \sigma\inv = ^{\k(\sigma)}\!u.
$$
By Skolem-Noether theorem, every automorphism of the
central simple $\Q$-algebra $\End^0(X)\cong\Mat_2(\H_p)$
is an inner one.
This implies that for each $\sigma\in G_1$ there exists
$w_\sigma \in\End^0(X)\mult$ such that
$$
\sigma u \sigma\inv = w_\sigma u w_\sigma\inv
$$
Since the center of $\End^0(X)$ is $\Q$, the choice of $w_\sigma$
is unique up to multiplication by a non-zero rational number. This implies
that $w_\sigma w_\tau$ equals $w_{\sigma\tau}$ times a non-zero rational
number.

Let 
$$
c'_\sigma=\sigma w_\sigma\inv.
$$
Each $c'_\sigma$ commutes with $\End^0(X)$ and hence with 
$\End^0(X) \qtimes \Q_2 = \End_{\Q_2}(V_2(X))$.
This means that $c'_\sigma \in \Q_2\mult \Id_{V_2(X)}$. The image
$$
c_\sigma \in \Q_2\mult  \Id_{V_2(X)} / \Q\mult  \Id_{V_2(X)} \cong \Q_2\mult / \Q\mult
$$
of $c'_\sigma$ in $\Q_2\mult / \Q\mult$ does not depend on the choice of $w_\sigma$.
Also, the map
$$
G_1 \to \Q_2\mult / \Q\mult, \sigma \mapsto c_\sigma
$$
is a group homomorphism. It has to be trivial, since $G_1$ is perfect. Therefore,
$$
c_\sigma \in \Q\mult \Id_{V_2(X)} \text{ for all }\sigma\in G_1.
$$
and
$$
\sigma = (c')\inv w_\sigma \in \End^0(X)\mult
$$
Finally, recall that $\End_{F_1}^0(X)=\i(D)$, so
$$
i(D) =
\{u\in\End^0(X)
\mid
^{\k(\sigma)}\!\!u=u\quad \forall \sigma \in \Gal(L_1/F_1)\},
$$
and each $\sigma \in G_1= \Gal(L_1/F_1)$ commutes with $\i(D)$.
\end{proof}

By combining $\rho: G_1 \hookrightarrow \End^0(X, \i)$ and $\i: D \hookrightarrow \End^0(X, \i)$
we obtain a natural homomorphism $D[G_1] \to \End^0(X, \i)$.

%%%%%%%%%%%%%%%%%%%%%%%%%%%%%%%%%%%%%%%%%%%%%%%%
\begin{claim}\label{cl5}
The $D$-algebra homomorphism $D[G_1] \to \End^0(X, \i)$ is surjective.
\end{claim}
\begin{proof}[Proof of Claim \ref{cl5}]
Let $M$ be the image of $D[G_1]$ in $\End^0(X, \i)$ under the above map.
Then $M \dtimes E$ coincides with the image of  $E[G_1] = D[G_1] \dtimes E$
in $\End^0(X, \i) \dtimes E =  \End_{E}(V_2(X))$. 
Since $E[G_1]$-module $V_2(X)$ is absolutely simple,
$$
E[G_1] \to \End_{E}(V_2(X))
$$
is surjective. Therefore,
$$
\dim_D(M) = \dim_D(\End^0(X, \i))
$$
and $M=\End^0(X, i)$, which proves the claim.
\end{proof}

%%%%%%%%%%%%%%%%%%%%%%%%%%%%%%%%%%%%%%%%%%%%%%%%
\begin{claim}\label{cl6}
The $D$-algebra $\End^0(X, \i)$ can be identified with a direct summand of $D[G_1]$.
\end{claim}
\begin{proof}[Proof of Claim \ref{cl6}]
The semisimplicity of $D[G_1]$ (by Maschke's Theorem)
and simplicity of $\End^0(X, \i)$ allow us to make such
an identification.
\end{proof}

This concludes the proof of Theorem \ref{nsauxthm}.
\end{proof}

%=========================================================================
\section{Examples in characteristic zero}
\label{examples}
%========================================================================= 
In this section we produce examples of 
hyperelliptic curves defined over $\Q$ with real 
quadratic fields $\Q(\sqrt{5})$ as endomorphism algebras 
of their jacobians. For clarity, we use capital 
Latin letters $B$, $C$, $D$, $T$ for indeterminates over 
a field $K$ and lower-case letters $b$, $c$, $d$, $t$ for their 
specializations in $K$. Further, we put $\eta=(\sqrt{5}-1)/2$,
so that $\Z[\eta]$ is the ring of integers of $\Q(\sqrt{5})$.

The following corollary of Theorem \ref{mainthm} will be used 
to obtain these examples.
\begin{cor}\label{maincor}
Let $K$ be a field of characteristic $0$ and 
$f_T(x)\in K(T)[x]$ be an irreducible separable polynomial
of degree $n=6$ in $x$ pa\-ra\-me\-trized by
a variable $T$ transcendental over $K$.
Define a hyperelliptic curve $C_T$ over $K(T)$ by
$$
C_T = C_{f_T}:y^2=f_T(x).
$$
Let $J(C_T)$ be its jacobian and
$\End_{K(T)}(J(C_T))$ be the ring of $K(T)$-endo\-morphisms of $J(C_T)$.
Assume that
\renewcommand{\theenumi}{\roman{enumi}}
\begin{enumerate}
\item $\Gal(f_T/K(T))\cong\A_5$,
\item $\End_{K(T)}^0(J(C_T))$ is isomorphic to 
a real quadratic field $D$, and
\item for some value $t\in K$ of $T$
the polynomial $f_t$ is irreducible 
over $K$ and $\Gal(f_t/K)\cong\A_5$.
\end{enumerate}
Then $\End(J(C_t))$ is isomorphic to an order of $D$ with
$$
\discr(\End(J(C_t)))\equiv5\bmod8.
$$
\end{cor}
\begin{proof}
Assume that all of the conditions of Corollary \ref{maincor}
are satisfied for specialization of $T$ to $t\in K$.
The action of $\Gal(f_T/K(T))$ on $\fR_T$ extends to
an action of $\Gal(f_T/K(T))$ on $\fR_t$.
The action can be factored through $\Gal(f_t/K)$;
and since $\Gal(f_T/K(T))\cong\Gal(f_t/K)$,
the $\Gal(f_T/K(T))$-sets $\fR_T$ and $\fR_t$ are also isomorphic.
We also have $\End_{K(T)}(J(C_T))\subset\End_K(J(C_t))$.
The conclusion of the corollary follows from Theorem \ref{mainthm}.
\end{proof}

In \cite{HashiOnBru}, K.~Hashimoto 
gives the following form of
Brumer's 3-parameter family of curves:
\begin{eqnarray}\label{brucurve}
\begin{aligned}
C_{B,C,D}: y^2 = f_{B,C,D}(x) = 
x^6 & +2Cx^5+(2+2C+C^2-4BD)x^4\\
&+(2+4B+2C+2C^2-4D-8BD)x^3\\
&+ (5+12B+4C+C^2-4BD)x^2\\
&+(6+12B+2C)x+4B+1.
\end{aligned}
\end{eqnarray}
For indeterminates $B, C, D$ over $\Q$,
the algebra of endomorphisms of 
its jacobian $J(C_{B,C,D})$ is isomorphic
to $\Q(\sqrt{5})$ and its endomorphisms ring
$\End(J(C_{B,C,D}))\cong\Z[\eta]$.
Moreover, the polynomial $f_{B,C,D}$ is irreducible
over $\Q(B,C,D)$, and its splitting field over $\Q(B,C,D)$ 
has $\A_5$ as its Galois group.

Assume that for certain values of $B,C,D$,
say $b, c, d\in\Q$, respectively, 
the polynomial $f_{b,c,d}$ is 
irreducible over $\Q$, and the Galois group 
of its splitting field over $\Q$ 
is isomorphic to $\A_5$.
Then $f_{B,c,d}, f_{B,C,d}$ are also
irreducible over $\Q(B)$ and $\Q(B,C)$,
respectively, because their factorization
would lead to a factorization of
$f_{b,c,d}$ over $\Q$. We
have a tower of groups
\begin{eqnarray*}
\A_5&\cong&\Gal(f_{b,c,d}/\Q)\\
&\subset&\Gal(f_{B,c,d}/\Q(B))\\
&\subset&\Gal(f_{B,C,d}/\Q(B,C))\\
&\subset&\Gal(f_{B,C,D}/\Q(B,C,D))\\
&\cong&\A_5,
\end{eqnarray*}
which forces every intermediate Galois 
group to be isomorphic to $\A_5$.
By applying Corollary \ref{maincor}
to the consecutive specializations
we conclude that every one of them has the 
ring of integers of the quadratic field 
$\Q(\sqrt{5})$ for the endomorphism ring.

By picking appropriate values $b, c, d$ 
we can give examples to Corollary \ref{maincor}. 
As an example, we work through 
one of the curves given in the table. 
The others can be treated in a similar 
manner.

%---------------------------------------------------------------------

\begin{eg}\label{eg1}
For $b=0$, $c=1$, $d=2$ we get a hyperelliptic curve
\begin{eqnarray*}
C: y^2 & = & f_{0, 1, 2}(x) \\
& = & {x^{6} + {2}x^{5} + {5}x^{4} - {2}x^{3} + {10}x^{2} + {8}x + 1}
\end{eqnarray*}
The splitting field of $f$ is an $\A_5$-extension of $\Q$
with discriminant $1125721=1061^2$.
By the above argument, we have $\End(J(C))\cong\Z[\eta]$.
By allowing one of the above variables to vary,
say $d=T$, while fixing the other two, we can get a 
hyperelliptic curve
\begin{eqnarray*}
C': y^2 & = & f_{0,1,T}(x) \\
& = & {x^6 + 2x^5 + 5x^4 + (6-4T)x^3 + 10x^2 + 8x + 1}
\end{eqnarray*}
defined over $\Q(T)$ such that $\Gal(f_T/\Q(T))\cong\A_5$ by the above argument. 
The endomorphism ring of its jacobian is also isomorphic to $\Z[\eta]$.
\end{eg}

Non-isogenuity of the jacobians of curves listed in Table \ref{tbl:egzero}
can be deduced from pairwise coprimality of the discriminants of the
splitting fields of polynomials that define them.

\begin{thm}
If polynomials $f(x), h(x)\in\Q[x]$ satisfy all of the conditions 
of Theorem \ref{mainthm}, and if the discriminants
of the splitting fields of $f$ and $h$ over $\Q$
are relatively prime, then the jacobians $J(C_f)$ and
$J(C_h)$ of the curves $C_f:y^2=f(x)$ and $C_h:y^2=h(x)$
are not isogenous.
\end{thm}
\begin{proof}
It is known that if $L/K$ is an extension of algebraic number
fields, then $\discr(K)|\discr(L)$ by the Discriminant Tower Theorem. 
This implies that the splitting fields $\Q(\fR_f)$ and $\Q(\fR_h)$
of given polynomials are linearly disjoint.
Indeed, if the two fields have a common subfield, say $E$, 
then $\discr(E)=1$ and $E=\Q$ by the theorem of Hermite.
By Theorem \ref{NICor}, $J(C_f)$ and $J(C_h)$ are not isogenous.
\end{proof}

By Hilbert's irreducibility theorem, there are infinitely many rational 
$t$ such that $\Gal(f_t/\Q)\cong\A_5$. Moreover, infinitely
many of these extensions are pairwise linearly disjoint over $\Q$. 
It follows from Theorem \ref{NICor} that there is an infinite number 
of pairwise non-isogenous hyperelliptic jacobians of dimension 2
with $\Z[\eta]$ as their endomorphism ring.

For each of the curves in Table \ref{tbl:egzero}, 
the polynomial $f$ is irreducible, 
and the splitting field of $f$ over $\Q$ is an $\A_5$-extension.
Note that the specialization
$b=1$, $c=1$, $d=2$ is the original example 
given in Brumer's paper \cite{BruRankJ}. This table was
constructed with the help of PARI-GP number-theoretic package \cite{PARI2}
through a search for integer values $-5\leq b, c, d\leq 5$.
The polynomials were chosen so that the discriminants of their
splitting fields are pairwise relatively prime in order to ensure that
the jacobian of the corresponding hyperelliptic curves are pairwise
nonisogenous.

More examples can be obtained by examination of polynomials 
in the table in Appendix A.3 of \cite{WilRM}.
J.~Wilson proved for every polynomial $f$ in that table, 
that if we define a hyperelliptic curve $C_f:y^2=f(x)$, then  
$\Z[\eta]\subset\End_{\Q}(J(C_f))$. Therefore, if $\Gal(f/\Q)\cong\A_5$, 
then $\End_{\Q}(J(C_f))=\End(J(C_f))=\Z[\eta]$.
For example, if $f(x)=3x^6+8x^5+54x^4-26x^3-173x^2+218x-73$, then
$\Gal(f/\Q)\cong\A_5$, and $\End(J(C_f))=\Z[\eta]$. It is easily
shown that this polynomial is not a specialization of the family 
(\ref{brucurve}) of Brumer's curve for any $b, c, d\in\Q$.
\vskip 1cm

\tablehead{
\hline
$b$ & $c$ & $d$ & $C:y^2=f_{b,c,d}(x)$ & $\discr(\Q(\fR_f))$ \\
\hline
\hline
}
\tabletail{\hline}
\topcaption{Some curves  \mbox{$C$} over  \mbox{$\Q$} with \mbox{$\End(J(C))=\Z[\eta]$}}
\begin{supertabular*}{\textwidth}
{|r@{\hspace*{0.2cm}}r@{\hspace*{0.2cm}}r@{\hspace*{0.3cm}\extracolsep{\fill}}l@{\hspace*{0.3cm}\extracolsep{\fill}}l|}\label{tbl:egzero}

$0$ & $0$ & $0$ &
$y^2={x^{6} + {2}x^{4} + {2}x^{3} + {5}x^{2} + {6}x + 1}$ &
$2^6\cdot103^2$\\

$0$ & $1$ & $2$ &
$y^2={x^{6} + {2}x^{5} + {5}x^{4} - {2}x^{3} + {10}x^{2} + {8}x + 1}$ &
$1061^2$\\

$0$ & $-1$ & $-3$ &
$y^2={x^{6} - {2}x^{5} + x^{4} + {14}x^{3} + {2}x^{2} + {4}x + 1}$ &
$11^2\cdot137^2$\\

$0$ & $-1$ & $5$ &
$y^2={x^{6} - {2}x^{5} + x^{4} - {18}x^{3} + {2}x^{2} + {4}x + 1}$ &
$2293^2$\\

$0$ & $3$ & $-3$ &
$y^2={x^{6} + {6}x^{5} + {17}x^{4} + {38}x^{3} + {26}x^{2} + {12}x + 1}$ &
$4483^2$\\

$0$ & $3$ & $5$ &
$y^2={x^{6} + {6}x^{5} + {17}x^{4} + {6}x^{3} + {26}x^{2} + {12}x + 1}$ &
$3^2\cdot4441^2$\\

$0$ & $-3$ & $-2$ &
$y^2={x^{6} - {6}x^{5} + {5}x^{4} + {22}x^{3} + {2}x^{2} + 1}$ &
$2609^2$\\

$0$ & $4$ & $2$ &
$y^2={x^{6} + {8}x^{5} + {26}x^{4} + {34}x^{3} + {37}x^{2} + {14}x + 1}$ &
$53^2\cdot79^2$\\

$0$ & $4$ & $-2$ &
$y^2={x^{6} + {8}x^{5} + {26}x^{4} + {50}x^{3} + {37}x^{2} + {14}x + 1}$ &
$2707^2$\\

$0$ & $-4$ & $2$ &
$y^2={x^{6} - {8}x^{5} + {10}x^{4} + {18}x^{3} + {5}x^{2} - {2}x + 1}$ &
$2029^2$\\

$0$ & $-4$ & $-2$ &
$y^2={x^{6} - {8}x^{5} + {10}x^{4} + {34}x^{3} + {5}x^{2} - {2}x + 1}$ &
$6827^2$\\

$0$ & $5$ & $2$ &
$y^2={x^{6} + {10}x^{5} + {37}x^{4} + {54}x^{3} + {50}x^{2} + {16}x + 1}$ &
$17^2\cdot337^2$\\

$0$ & $5$ & $-2$ &
$y^2={x^{6} + {10}x^{5} + {37}x^{4} + {70}x^{3} + {50}x^{2} + {16}x + 1}$ &
$5^2\cdot757^2$\\

$0$ & $-5$ & $1$ &
$y^2={x^{6} - {10}x^{5} + {17}x^{4} + {38}x^{3} + {10}x^{2} - {4}x + 1}$ &
$3929^2$\\

$0$ & $-5$ & $5$ &
$y^2={x^{6} - {10}x^{5} + {17}x^{4} + {22}x^{3} + {10}x^{2} - {4}x + 1}$ &
$47^2\cdot251^2$\\

$1$ & $0$ & $4$ &
$y^2={x^{6} - {14}x^{4} - {42}x^{3} + x^{2} + {18}x + {5}}$ &
$41^2\cdot941^2$\\

$1$ & $-1$ & $2$ &
$y^2={x^{6} - {2}x^{5} - {7}x^{4} - {18}x^{3} + {6}x^{2} + {16}x + {5}}$ &
$7933^2$\\

$1$ & $-1$ & $-2$ &
$y^2={x^{6} - {2}x^{5} + {9}x^{4} + {30}x^{3} + {22}x^{2} + {16}x + {5}}$ &
$19^2\cdot1289^2$\\

$1$ & $2$ & $2$ &
$y^2={x^{6} + {4}x^{5} + {2}x^{4} - {6}x^{3} + {21}x^{2} + {22}x + {5}}$ &
$2861^2$\\

$1$ & $-2$ & $2$ &
$y^2={x^{6} - {4}x^{5} - {6}x^{4} - {14}x^{3} + {5}x^{2} + {14}x + {5}}$ &
$9907^2$\\

$1$ & $-2$ & $-2$ &
$y^2={x^{6} - {4}x^{5} + {10}x^{4} + {34}x^{3} + {21}x^{2} + {14}x + {5}}$ &
$71^2\cdot607^2$\\

$1$ & $-3$ & $0$ &
$y^2={x^{6} - {6}x^{5} + {5}x^{4} + {18}x^{3} + {14}x^{2} + {12}x + {5}}$ &
$3089^2$\\

$1$ & $-3$ & $4$ &
$y^2={x^{6} - {6}x^{5} - {11}x^{4} - {30}x^{3} - {2}x^{2} + {12}x + {5}}$ &
$23^2\cdot3137^2$\\

$1$ & $4$ & $4$ &
$y^2={x^{6} + {8}x^{5} + {10}x^{4} - {2}x^{3} + {33}x^{2} + {26}x + {5}}$ &
$17509^2$\\

$1$ & $4$ & $-4$ &
$y^2={x^{6} + {8}x^{5} + {42}x^{4} + {94}x^{3} + {65}x^{2} + {26}x + {5}}$ &
$55763^2$\\

$1$ & $5$ & $-4$ &
$y^2={x^{6} + {10}x^{5} + {53}x^{4} + {114}x^{3} + {78}x^{2} + {28}x + {5}}$ &
$55793^2$\\

$1$ & $-5$ & $2$ &
$y^2={x^{6} - {10}x^{5} + {9}x^{4} + {22}x^{3} + {14}x^{2} + {8}x + {5}}$ &
$34729^2$\\

$-1$ & $1$ & $3$ &
$y^2={x^{6} + {2}x^{5} + {17}x^{4} + {14}x^{3} + {10}x^{2} - {4}x - {3}}$ &
$11027^2$\\

$-1$ & $1$ & $-5$ &
$y^2={x^{6} + {2}x^{5} - {15}x^{4} - {18}x^{3} - {22}x^{2} - {4}x - {3}}$ &
$17387^2$\\

$-1$ & $-1$ & $3$ &
$y^2={x^{6} - {2}x^{5} + {13}x^{4} + {10}x^{3} + {2}x^{2} - {8}x - {3}}$ &
$9293^2$\\

$-1$ & $-1$ & $-5$ &
$y^2={x^{6} - {2}x^{5} - {19}x^{4} - {22}x^{3} - {30}x^{2} - {8}x - {3}}$ &
$26501^2$\\

$-1$ & $-2$ & $3$ &
$y^2={x^{6} - {4}x^{5} + {14}x^{4} + {14}x^{3} + x^{2} - {10}x - {3}}$ &
$13^2\cdot1151^2$\\

$-1$ & $4$ & $3$ &
$y^2={x^{6} + {8}x^{5} + {38}x^{4} + {50}x^{3} + {37}x^{2} + {2}x - {3}}$ &
$157^2\cdot389^2$\\

$-1$ & $4$ & $-5$ &
$y^2={x^{6} + {8}x^{5} + {6}x^{4} + {18}x^{3} + {5}x^{2} + {2}x - {3}}$ &
$43^2\cdot227^2$\\

$-1$ & $-4$ & $3$ &
$y^2={x^{6} - {8}x^{5} + {22}x^{4} + {34}x^{3} + {5}x^{2} - {14}x - {3}}$ &
$59^2\cdot1483^2$\\

$-1$ & $-4$ & $-5$ &
$y^2={x^{6} - {8}x^{5} - {10}x^{4} + {2}x^{3} - {27}x^{2} - {14}x - {3}}$ &
$83417^2$\\

$-1$ & $5$ & $3$ &
$y^2={x^{6} + {10}x^{5} + {49}x^{4} + {70}x^{3} + {50}x^{2} + {4}x - {3}}$ &
$73^2\cdot1511^2$\\

$-1$ & $5$ & $-5$ &
$y^2={x^{6} + {10}x^{5} + {17}x^{4} + {38}x^{3} + {18}x^{2} + {4}x - {3}}$ &
$50423^2$\\

$2$ & $0$ & $-3$ &
$y^2={x^{6} + {26}x^{4} + {70}x^{3} + {53}x^{2} + {30}x + {9}}$ &
$167^2\cdot181^2$\\

$2$ & $-1$ & $4$ &
$y^2={x^{6} - {2}x^{5} - {31}x^{4} - {70}x^{3} - {6}x^{2} + {28}x + {9}}$ &
$61^2\cdot2593^2$\\

$2$ & $-1$ & $-4$ &
$y^2={x^{6} - {2}x^{5} + {33}x^{4} + {90}x^{3} + {58}x^{2} + {28}x + {9}}$ &
$457^2\cdot2011^2$\\

$2$ & $-2$ & $4$ &
$y^2={x^{6} - {4}x^{5} - {30}x^{4} - {66}x^{3} - {7}x^{2} + {26}x + {9}}$ &
$107^2\cdot2693^2$\\

$2$ & $-2$ & $-4$ &
$y^2={x^{6} - {4}x^{5} + {34}x^{4} + {94}x^{3} + {57}x^{2} + {26}x + {9}}$ &
$1280761^2$\\

$2$ & $4$ & $5$ &
$y^2={x^{6} + {8}x^{5} - {14}x^{4} - {50}x^{3} + {21}x^{2} + {38}x + {9}}$ &
$80407^2$\\

$2$ & $-4$ & $-3$ &
$y^2={x^{6} - {8}x^{5} + {34}x^{4} + {94}x^{3} + {53}x^{2} + {22}x + {9}}$ &
$673^2\cdot2087^2$\\

$2$ & $5$ & $-3$ &
$y^2={x^{6} + {10}x^{5} + {61}x^{4} + {130}x^{3} + {98}x^{2} + {40}x + {9}}$ &
$31^4\cdot233^2$\\

$2$ & $5$ & $5$ &
$y^2={x^{6} + {10}x^{5} - {3}x^{4} - {30}x^{3} + {34}x^{2} + {40}x + {9}}$ &
$145487^2$\\

$2$ & $-5$ & $0$ &
$y^2={x^{6} - {10}x^{5} + {17}x^{4} + {50}x^{3} + {34}x^{2} + {20}x + {9}}$ &
$29663^2$\\

$2$ & $-5$ & $4$ &
$y^2={x^{6} - {10}x^{5} - {15}x^{4} - {30}x^{3} + {2}x^{2} + {20}x + {9}}$ &
$271273^2$\\

$-2$ & $0$ & $-3$ &
$y^2={x^{6} - {22}x^{4} - {42}x^{3} - {43}x^{2} - {18}x - {7}}$ &
$32771^2$\\

$-2$ & $-1$ & $-4$ &
$y^2={x^{6} - {2}x^{5} - {31}x^{4} - {54}x^{3} - {54}x^{2} - {20}x - {7}}$ &
$29^2\cdot3319^2$\\

$-2$ & $-2$ & $4$ &
$y^2={x^{6} - {4}x^{5} + {34}x^{4} + {46}x^{3} + {9}x^{2} - {22}x - {7}}$ &
$11057^2$\\

$-2$ & $3$ & $4$ &
$y^2={x^{6} + {6}x^{5} + {49}x^{4} + {66}x^{3} + {34}x^{2} - {12}x - {7}}$ &
$582983^2$\\

$-2$ & $-3$ & $-3$ &
$y^2={x^{6} - {6}x^{5} - {19}x^{4} - {30}x^{3} - {46}x^{2} - {24}x - {7}}$ &
$72101^2$\\

$-2$ & $4$ & $5$ &
$y^2={x^{6} + {8}x^{5} + {66}x^{4} + {94}x^{3} + {53}x^{2} - {10}x - {7}}$ &
$7^2\cdot228281^2$\\

$-2$ & $-4$ & $-3$ &
$y^2={x^{6} - {8}x^{5} - {14}x^{4} - {18}x^{3} - {43}x^{2} - {26}x - {7}}$ &
$144223^2$\\

$-2$ & $5$ & $5$ &
$y^2={x^{6} + {10}x^{5} + {77}x^{4} + {114}x^{3} + {66}x^{2} - {8}x - {7}}$ &
$317^2\cdot7057^2$\\

$3$ & $0$ & $-5$ &
$y^2={x^{6} + {62}x^{4} + {154}x^{3} + {101}x^{2} + {42}x + {13}}$ &
$5562929^2$\\

$3$ & $1$ & $3$ &
$y^2={x^{6} + {2}x^{5} - {31}x^{4} - {66}x^{3} + {10}x^{2} + {44}x + {13}}$ &
$67^2\cdot2617^2$\\

$3$ & $-1$ & $-5$ &
$y^2={x^{6} - {2}x^{5} + {61}x^{4} + {154}x^{3} + {98}x^{2} + {40}x + {13}}$ &
$6683357^2$\\

$3$ & $2$ & $-5$ &
$y^2={x^{6} + {4}x^{5} + {70}x^{4} + {166}x^{3} + {113}x^{2} + {46}x + {13}}$ &
$223^2\cdot20549^2$\\

$3$ & $-3$ & $3$ &
$y^2={x^{6} - {6}x^{5} - {31}x^{4} - {58}x^{3} + {2}x^{2} + {36}x + {13}}$ &
$384641^2$\\

$3$ & $-4$ & $-5$ &
$y^2={x^{6} - {8}x^{5} + {70}x^{4} + {178}x^{3} + {101}x^{2} + {34}x + {13}}$ &
$14931629^2$\\

$3$ & $-5$ & $-1$ &
$y^2={x^{6} - {10}x^{5} + {29}x^{4} + {82}x^{3} + {58}x^{2} + {32}x + {13}}$ &
$229^2\cdot2897^2$\\

$-3$ & $1$ & $-4$ &
$y^2={x^{6} + {2}x^{5} - {43}x^{4} - {86}x^{3} - {74}x^{2} - {28}x - {11}}$ &
$269^2\cdot3109^2$\\

$-3$ & $-3$ & $4$ &
$y^2={x^{6} - {6}x^{5} + {53}x^{4} + {82}x^{3} + {14}x^{2} - {36}x - {11}}$ &
$1493^2\cdot2161^2$\\

$-3$ & $-4$ & $-4$ &
$y^2={x^{6} - {8}x^{5} - {38}x^{4} - {66}x^{3} - {79}x^{2} - {38}x - {11}}$ &
$238991^2$\\

$-3$ & $5$ & $-4$ &
$y^2={x^{6} + {10}x^{5} - {11}x^{4} - {30}x^{3} - {34}x^{2} - {20}x - {11}}$ &
$236813^2$\\

$-3$ & $-5$ & $2$ &
$y^2={x^{6} - {10}x^{5} + {41}x^{4} + {70}x^{3} - {2}x^{2} - {40}x - {11}}$ &
$1066643^2$\\

$4$ & $-1$ & $5$ &
$y^2={x^{6} - {2}x^{5} - {79}x^{4} - {162}x^{3} - {30}x^{2} + {52}x + {17}}$ &
$991^2\cdot5441^2$\\

$4$ & $-3$ & $-2$ &
$y^2={x^{6} - {6}x^{5} + {37}x^{4} + {102}x^{3} + {82}x^{2} + {48}x + {17}}$ &
$1776617^2$\\

$4$ & $-4$ & $2$ &
$y^2={x^{6} - {8}x^{5} - {22}x^{4} - {30}x^{3} + {21}x^{2} + {46}x + {17}}$ &
$37^2\cdot8161^2$\\

$4$ & $-4$ & $-2$ &
$y^2={x^{6} - {8}x^{5} + {42}x^{4} + {114}x^{3} + {85}x^{2} + {46}x + {17}}$ &
$331^2\cdot8297^2$\\

$4$ & $-5$ & $1$ &
$y^2={x^{6} - {10}x^{5} + x^{4} + {22}x^{3} + {42}x^{2} + {44}x + {17}}$ &
$123829^2$\\

$-4$ & $-1$ & $5$ &
$y^2={x^{6} - {2}x^{5} + {81}x^{4} + {126}x^{3} + {34}x^{2} - {44}x - {15}}$ &
$89^2\cdot154127^2$\\

$-4$ & $2$ & $5$ &
$y^2={x^{6} + {4}x^{5} + {90}x^{4} + {138}x^{3} + {49}x^{2} - {38}x - {15}}$ &
$139^2\cdot94651^2$\\

$-4$ & $-5$ & $-3$ &
$y^2={x^{6} - {10}x^{5} - {31}x^{4} - {58}x^{3} - {86}x^{2} - {52}x - {15}}$ &
$500111^2$\\

$5$ & $-4$ & $0$ &
$y^2={x^{6} - {8}x^{5} + {10}x^{4} + {46}x^{3} + {65}x^{2} + {58}x + {21}}$ &
$65003^2$\\

$-5$ & $1$ & $-5$ &
$y^2={x^{6} + {2}x^{5} - {95}x^{4} - {194}x^{3} - {150}x^{2} - {52}x - {19}}$ &
$3301^2\cdot5783^2$\\

$-5$ & $-1$ & $-5$ &
$y^2={x^{6} - {2}x^{5} - {99}x^{4} - {198}x^{3} - {158}x^{2} - {56}x - {19}}$ &
$15327437^2$\\

$-5$ & $-4$ & $3$ &
$y^2={x^{6} - {8}x^{5} + {70}x^{4} + {114}x^{3} + {5}x^{2} - {62}x - {19}}$ &
$9647317^2$\\

%$-5$ & $-5$ & $-1$ &
%$y^2={x^{6} - {10}x^{5} - {3}x^{4} - {14}x^{3} - {70}x^{2} - {64}x - {19}}$ &
%$401477^2$\\
\end{supertabular*}

%====================================================================
\section{Examples in Positive Characteristic}\label{PosEg}
%====================================================================
Examples of hyperelliptic curves 
satisfying conditions of Theorem \ref{mainthm} can be found by reducing examples
defined over $\Q(T)$ (for example,  specializations of the Brumer's family) 
modulo odd primes. We need to make sure that
after reduction the defining polynomial remains irreducible, separable,
and has Galois group $\A_5$.

%=================== Testing for supersingularity =======================

In positive characteristic it becomes necessary to distinguish between the
two outcomes of that theorem. Let $\chr(K)=p>2$,
$$
f(x) = x^6+a_5 x^5+a_4 x^4+a_3 x^3+a_2 x^2+a_1 x + a_0,
$$
be a polynomial with distinct roots over $K_a$, and let $f(x)^{(p-1)/2}=\sum c_i x^i$. 
Then the Cartier-Manin/Hasse-Witt matrix \cite{NygCC, YuiJV} for the hyperelliptic curve 
$C: y^2=f(x)$ is the matrix obtained from
$$
M=
\begin{pmatrix}
c_{p-1} & c_{p-2} \\
c_{2p-1} & c_{2p-2}
\end{pmatrix}.
$$
by extraction of $p$th roots of the entries.
It is known \cite{ManCFG}, \cite[Th.~2.2]{YuiJV}, \cite{IbuSC} that the jacobian of the curve $C$ is a supersingular
abelian variety if and only if
$$
M M^{(p)}=
\begin{pmatrix}
c_{p-1} & c_{p-2} \\
c_{2p-1} & c_{2p-2}
\end{pmatrix}
\cdot
\begin{pmatrix}
c_{p-1}^p & c_{p-2}^p \\
c_{2p-1}^p & c_{2p-2}^p
\end{pmatrix}
=0.
$$
It can be seen that this happens in one of the following cases:
\begin{enumerate}
\item $c_{p-1}=c_{2p-1}=c_{2p-2}=0$, or
\item $c_{2p-1}\neq 0$, 
$$
c_{p-2}=-\frac{c_{p-1}^{p+1}}{c_{2p-1}^{p}}
\quad\text{and}\quad
c_{2p-2}=-\frac{c_{p-1}^{p}}{c_{2p-1}^{p-1}}.
$$
\end{enumerate}

%================================================================

\begin{eg}\label{charthree}
Let $p=\chr(K)=3$. Then $f(x)^{(p-1)/2}=f(x)$, 
and the curve $C: y^2=f(x)$ is supersingular if and only if one of the following conditions holds:
\begin{enumerate}
\item $a_2=a_4=a_5=0$. In this case, $a_1\neq 0$, since polynomials
$f(x)=x^6+a_3 x^3 + a_0$ are not separable in characteristic $3$.
\item $a_5\neq 0$, $a_1=-a_2^4/a_5^3$, and
$a_4=-a_2^3/a_5^2$.
\end{enumerate}

This means $C$ is supersingular if and only if either
$$
f(x) = x^6 + a_3 x^3 +  a_1 x + a_0,
\quad a_1\neq 0,
$$
or
$$
f(x) = x^6 + a_5 x^5 - \frac{a_2^3}{a_5^2} x^4 + a_3 x^3 + a_2 x^2 - \frac{a_2^4}{a_5^3} x + a_0, 
\quad a_5\neq0.
$$
\end{eg}

\begin{eg}\label{morecharthree}
Let $C_f$ be the smooth hyperelliptic curve $y^2=f(x)$ over $\F_3(T)$,
where $f(x)$ is one of the polynomials in Table \ref{tbl:eg3nonss}. Then 
$$
\End(J(C_f)) = \End_{\F_3(T)}(J(C_f)) = \Z[\eta].
$$
\vskip 0.5cm

\begin{center}
\tablehead{
\hline
$b$ & $c$ & $d$ & $C:y^2=f_{b,c,d}(x)$ reduced modulo $3$ \\
\hline
\hline
}
\tabletail{\hline}
\topcaption{Some curves  \mbox{$C$} over  \mbox{$\F_3(T)$} with \mbox{$\End(J(C))=\Z[\eta]$}}
\begin{supertabular}
{|r@{\hspace*{0.5cm}}r@{\hspace*{0.5cm}}r@{\hspace*{0.5cm}\extracolsep{\fill}}l|}
\label{tbl:eg3nonss}
$0$ & $1$ & $T$ & $y^2 = x^6 + 2x^5 + 2x^4 + 2Tx^3 + x^2 + 2x + 1$\\
$0$ & $2$ & $T + 2$ & $y^2 = x^6 + x^5 + x^4 + 2Tx^3 + 2x^2 + x + 1$\\
$1$ & $1$ & $2T + 1$ & $y^2 = x^6 + 2x^5 + (T + 1)x^4 + x^3 + Tx^2 + 2x + 2$\\
%$1$ & $1$ & $2T + 2$ & $y^2 = x^6 + 2x^5 + Tx^4 + x^3 + (T + 2)x^2 + 2x + 2$\\
$1$ & $2$ & $T + 1$ & $y^2 = x^6 + x^5 + 2Tx^4 + (2T + 1)x^2 + x + 2$\\
%$1$ & $2$ & $2T + 1$ & $y^2 = x^6 + x^5 + Tx^4 + (T + 1)x^2 + x + 2$\\
$T$ & $1$ & $0$ & $y^2 = x^6 + 2x^5 + 2x^4 + Tx^3 + x^2 + 2x + T + 1$\\
%$T$ & $2$ & $0$ & $y^2 = x^6 + x^5 + x^4 + (T + 2)x^3 + 2x^2 + x + T + 1$\\
$T + 1$ & $0$ & $T + 2$ & $y^2 = x^6 + 2T^2x^4 + T^2x^3 + 2T^2x^2 + T + 2$\\
%$T + 1$ & $2$ & $0$ & $y^2 = x^6 + x^5 + x^4 + Tx^3 + 2x^2 + x + T + 2$\\
$T + 1$ & $2T$ & $T + 1$ & $y^2 = x^6 + Tx^5 + (2T + 1)x^4 + x^2 + Tx + T + 2$\\
%$T + 2$ & $1$ & $0$ & $y^2 = x^6 + 2x^5 + 2x^4 + (T + 2)x^3 + x^2 + 2x + T$\\
$T + 2$ & $2$ & $0$ & $y^2 = x^6 + x^5 + x^4 + (T + 1)x^3 + 2x^2 + x + T$\\
$T + 2$ & $2T + 2$ & $T + 2$ & $y^2 = x^6 + (T + 1)x^5 + 2Tx^4 + x^2 + (T + 
1)x + T$\\
$2T + 1$ & $0$ & $2T + 2$ & $y^2 = x^6 + 2T^2x^4 + T^2x^3 + 2T^2x^2 + 2T
+ 2$\\
%$2T + 1$ & $2$ & $0$ & $y^2 = x^6 + x^5 + x^4 + 2Tx^3 + 2x^2 + x + 2T + 2$\\
\end{supertabular}
\end{center}
\vskip 0.5cm

Let us work through one of the examples. 
The smooth curve $C$ defined over $\F_3(T)$ by
$$
C: y^2 = x^6 + 2x^5 + 2x^4 + 2Tx^3 + x^2 + 2x + 1
$$
is the reduction modulo $3$ of the curve
$$
C' = C_{0, 1, T}: y^2 = x^6 + 2x^5 + 5x^4 + (6-4T)x^3 + 10x^2 + 8x + 1
$$
defined over $\Q(T)$ which, according to Example \ref{eg1}, satisfies 
$$
\End(J(C')) = \End_{\Q(T)}(J(C')) \cong \Z[\eta].
$$
Thus we have
$$
\Z[\eta] \cong \End_{\Q(T)}(J(C')) \hookrightarrow \End_{\F_3(T)}(J(C)).
$$
Using MAGMA Computational Algebra System \cite{MAGMA, StPM}
we verify that the polynomial $f(x)= x^6 + 2x^5 + 2x^4 + 2Tx^3 + x^2 + 2x + 1$
is irreducible and separable over $\F_3(T)$ with $\Gal(f/\F_3(T))\cong \A_5$.
Finally, the procedure delineated above shows that 
$J(C)$ is not a supersingular abelian variety.
Therefore, by Theorem \ref{mainthm} we have $\End(J(C)) \cong \Z[\eta]$.

Other  examples from Table \ref{tbl:eg3nonss}, as well as Tables
\ref{tbl:eg5nonss},  \ref{tbl:eg7nonss}, and \ref{tbl:eg11nonss},
are obtained in a similar fashion.
Table \ref{tbl:eg5ss} is generated analogously, with the distinction that 
only supersingular examples are selected.
\end{eg}
%------------------------------------
All jacobians of smooth hyperelliptic curves over $\F_3(T)$ obtained by 
reduction of specializations of Brumer's family will be non-supersingular.
To prove this, observe that the reduction of the Brumer's equation \ref{brucurve}
modulo $3$ yields
$$
a_1=a_5\quad \text{ and }\quad a_4=a_2+C = a_2 - a_1
$$
This immediately rules out the first case of supersingularity
outlined in Example \ref{charthree}. 
In the second case, $a_1=a_5\neq 0$ and
$$
a_1 a_5^3= - a_2^4,
$$
so
$$
a_1^4= - a_2^4.
$$
For $a_1, a_2\in \F_3(T)$, this equation does not have a solution.

It is possible that supersingular examples exist over $\F_9(T)$.
Indeed, $a_2 = \varepsilon a_1$, where $\varepsilon$ is a root of $z^4+1=0$ in $\overline\F_3$.
We also know that $a_4 a_5^2 = - a_2^3$, so $a_4 =  -\varepsilon^3 a_1$.
Finally, plug this into $a_4 = a_2 - a_1$ and divide by $a_1\neq 0$, 
to get $\varepsilon^3+\varepsilon-1=0$.
The simultaneous solutions of $z^4+1=0$ and $z^3+z-1=0$ are
$$
\varepsilon = -1\pm\sqrt{-1} \in \F_9.
$$
Since $a_4 = 2+2C+C^2-4BD$, $a_5=a_1=2C$, and $a_4=(\varepsilon-1) a_1$,
we have
$$
BD = C^2+(\varepsilon+1)C-1.
$$
As a result, the equation of the curve will have the form
\begin{align*}
C: y^2  =  f(x) =  x^6 & - Cx^5 + (1-\varepsilon) Cx^4 \\
 & + (1+\varepsilon C+B-D)x^3 - \varepsilon C x^2 - Cx +(B+1).
\end{align*}
In order to verify that $\End_{\F_9(T)}(J(C)) = \Z[\eta]$, one will need to check that
$C$ is a reduction of a curve satisfying the conditions of Theorem \ref{mainthm},
that $f(x)$ is irreducible in $\F_9(T)[x]$, and that $\Gal(f/\F_9(T))\cong \A_5$.

When $\chr(K)>3$, any irreducible polynomial $f(x)$ of degree $6$ is separable, and
therefore the curve $y^2=f(x)$ defined by such a polynomial is smooth.

\begin{eg}
Let $C_f$ be the smooth hyperelliptic curve $y^2=f(x)$ over $\F_5(T)$,
where $f(x)$ is one of the polynomials in Table \ref{tbl:eg5nonss}. Then 
$$
\End(J(C_f)) = \End_{\F_5(T)}(J(C_f)) = \Z[\eta].
$$
\begin{center}
\tablehead{
\hline
$b$ & $c$ & $d$ & $C:y^2=f_{b,c,d}(x)$ reduced modulo $5$ \\
\hline
\hline
}
\tabletail{\hline}
\topcaption{Some curves  \mbox{$C$} over  \mbox{$\F_5(T)$} with \mbox{$\End(J(C))=\Z[\eta]$}}
\begin{supertabular}
{|r@{\hspace*{0.5cm}}r@{\hspace*{0.5cm}}r@{\hspace*{0.5cm}\extracolsep{\fill}}l|}
\label{tbl:eg5nonss}
$0$ & $0$ & $T + 3$ & $y^2 = x^6 + 2x^4 + Tx^3 + x + 1$\\
$0$ & $2$ & $T$ & $y^2 = x^6 + 4x^5 + (T + 4)x^3 + 2x^2 + 1$\\
$0$ & $3$ & $T + 4$ & $y^2 = x^6 + x^5 + 2x^4 + Tx^3 + x^2 + 2x + 1$\\
$0$ & $4$ & $T + 3$ & $y^2 = x^6 + 3x^5 + x^4 + Tx^3 + 2x^2 + 4x + 1$\\
$3$ & $3$ & $2T + 1$ & $y^2 = x^6 + x^5 + Tx^4 + 4Tx^3 + Tx^2 + 3x + 3$\\
$4$ & $0$ & $4T + 3$ & $y^2 = x^6 + (T + 4)x^4 + Tx^3 + Tx^2 + 4x + 2$\\
$4$ & $2$ & $4T$ & $y^2 = x^6 + 4x^5 + Tx^4 + Tx^3 + Tx^2 + 3x + 2$\\
$4$ & $3$ & $4T + 2$ & $y^2 = x^6 + x^5 + Tx^4 + Tx^3 + (T + 2)x^2 + 2$\\
%$T$ & $0$ & $3$ & $y^2 = x^6 + (3T + 2)x^4 + (2T + 1)x + 4T + 1$\\
%$T + 1$ & $0$ & $3$ & $y^2 = x^6 + 3Tx^4 + (2T + 3)x + 4T$\\
%$T + 2$ & $0$ & $3$ & $y^2 = x^6 + (3T + 3)x^4 + 2Tx + 4T + 4$\\
$T + 2$ & $4T$ & $4T + 4$ & $y^2 = x^6 + 3Tx^5 + 2x^2 + 4T + 4$\\
$T + 3$ & $4$ & $3$ & $y^2 = x^6 + 3x^5 + 3Tx^4 + 2x^2 + 2Tx + 4T + 3$\\
%$T + 4$ & $0$ & $3$ & $y^2 = x^6 + (3T + 4)x^4 + (2T + 4)x + 4T + 2$\\
%$2T$ & $0$ & $3$ & $y^2 = x^6 + (T + 2)x^4 + (4T + 1)x + 3T + 1$\\
%$2T + 1$ & $0$ & $3$ & $y^2 = x^6 + Tx^4 + (4T + 3)x + 3T$\\
%$2T + 2$ & $0$ & $3$ & $y^2 = x^6 + (T + 3)x^4 + 4Tx + 3T + 4$\\
$2T + 2$ & $3T$ & $3T + 4$ & $y^2 = x^6 + Tx^5 + 2x^2 + 3T + 4$\\
%$2T + 3$ & $4$ & $3$ & $y^2 = x^6 + 3x^5 + Tx^4 + 2x^2 + 4Tx + 3T + 3$\\
%$2T + 4$ & $0$ & $3$ & $y^2 = x^6 + (T + 4)x^4 + (4T + 4)x + 3T + 2$\\
%$3T$ & $0$ & $3$ & $y^2 = x^6 + (4T + 2)x^4 + (T + 1)x + 2T + 1$\\
%$3T + 1$ & $0$ & $3$ & $y^2 = x^6 + 4Tx^4 + (T + 3)x + 2T$\\
$3T + 1$ & $2T + 4$ & $2T + 1$ & $y^2 = x^6 + (4T + 3)x^5 + 2x^4 + 4x^3 +
x + 2T$\\
$3T + 2$ & $0$ & $3$ & $y^2 = x^6 + (4T + 3)x^4 + Tx + 2T + 4$\\
$3T + 2$ & $2T$ & $2T + 4$ & $y^2 = x^6 + 4Tx^5 + 2x^2 + 2T + 4$\\
%$3T + 3$ & $4$ & $3$ & $y^2 = x^6 + 3x^5 + 4Tx^4 + 2x^2 + Tx + 2T + 3$\\
$3T + 3$ & $2T + 2$ & $2T + 4$ & $y^2 = x^6 + (4T + 4)x^5 + 2x^4 + 4x^3 +
x + 2T + 3$\\
%$3T + 3$ & $2T + 4$ & $2T + 3$ & $y^2 = x^6 + (4T + 3)x^5 + 2x^2 + 2T + 3$\\
%$3T + 4$ & $0$ & $3$ & $y^2 = x^6 + (4T + 4)x^4 + (T + 4)x + 2T + 2$\\
$3T + 4$ & $2T + 1$ & $2T + 3$ & $y^2 = x^6 + (4T + 2)x^5 + 2x^4 + 4x^3 +
x + 2T + 2$\\
$3T + 4$ & $2T + 3$ & $2T + 2$ & $y^2 = x^6 + (4T + 1)x^5 + 2x^2 + 2T + 
2$\\
%$4T$ & $0$ & $3$ & $y^2 = x^6 + (2T + 2)x^4 + (3T + 1)x + T + 1$\\
%$4T + 1$ & $0$ & $3$ & $y^2 = x^6 + 2Tx^4 + (3T + 3)x + T$\\
$4T + 1$ & $1$ & $3$ & $y^2 = x^6 + 2x^5 + (2T + 3)x^4 + 4x^3 + 3Tx + 
T$\\
$4T + 1$ & $T + 4$ & $T + 1$ & $y^2 = x^6 + (2T + 3)x^5 + 2x^4 + 4x^3 + x +
T$\\
$4T + 1$ & $2T$ & $4T + 3$ & $y^2 = x^6 + 4Tx^5 + Tx^3 + 2Tx^2 + (2T + 3)x + T$\\
%$4T + 2$ & $0$ & $3$ & $y^2 = x^6 + (2T + 3)x^4 + 3Tx + T + 4$\\
$4T + 2$ & $T$ & $T + 4$ & $y^2 = x^6 + 2Tx^5 + 2x^2 + T + 4$\\
%$4T + 3$ & $4$ & $3$ & $y^2 = x^6 + 3x^5 + 2Tx^4 + 2x^2 + 3Tx + T + 3$\\
$4T + 3$ & $T + 2$ & $T + 4$ & $y^2 = x^6 + (2T + 4)x^5 + 2x^4 + 4x^3 + x + T + 3$\\
%$4T + 3$ & $T + 4$ & $T + 3$ & $y^2 = x^6 + (2T + 3)x^5 + 2x^2 + T + 3$\\
\end{supertabular}
\end{center}
\vskip 0.5cm
\end{eg}

\begin{eg}
Let $C_f$ be the smooth hyperelliptic curve $y^2=f(x)$ over $\F_5(T)$,
where $f(x)$ is one of the polynomials in Table \ref{tbl:eg5ss}. Then 
$$
\End_{\F_5(T)}(J(C_f)) = \Z[\eta].
$$
but $J(C)$ is a supersingular abelian variety; that is, $\End^0(J(C_f)) = \Mat_2(\H_5)$.
\begin{center}
\tablehead{
\hline
$b$ & $c$ & $d$ & $C:y^2=f_{b,c,d}(x)$ reduced modulo $5$ \\
\hline
\hline
}
\tabletail{\hline}
\topcaption{Some supersingular curves  \mbox{$C$} over  \mbox{$\F_5(T)$} 
with \mbox{$\End_{\F_5(T)}(J(C))=\Z[\eta]$}}
\begin{supertabular}
{|r@{\hspace*{0.3cm}}r@{\hspace*{0.3cm}}r@{\hspace*{0.3cm}\extracolsep{\fill}}l|}
\label{tbl:eg5ss}
$0$ & $1$ & $T$ & $y^2 = x^6 + 2x^5 + (T + 1)x^3 + 3x + 1$\\
%$T$ & $2$ & $0$ & $y^2 = x^6 + 4x^5 + (4T + 4)x^3 + (2T + 2)x^2 + 2Tx + 4T + 1$\\
$T$ & $4T + 3$ & $4T + 3$ & $y^2 = x^6 + (3T + 1)x^5 + 2x^4 + 4x^3 + x^2 + 2x + 4T + 1$\\
%$T + 1$ & $2$ & $0$ & $y^2 = x^6 + 4x^5 + (4T + 3)x^3 + (2T + 4)x^2 + (2T + 2)x + 4T$\\
%$T + 1$ & $T$ & $4T + 3$ & $y^2 = x^6 + 2Tx^5 + 4Tx^4 + 4Tx^3 + 3Tx^2 + (4T + 3)x + 4T$\\
%$T + 1$ & $T + 1$ & $4T + 3$ & $y^2 = x^6 + (2T + 2)x^5 + (T + 3)x^4 + (3T + 4)x^3 + 4Tx + 4T$\\
$T + 1$ & $4T + 2$ & $4T + 2$ & $y^2 = x^6 + (3T + 4)x^5 + 2x^4 + 4x^3 + x^2 + 2x + 4T$\\
$T + 2$ & $4T + 1$ & $4T + 1$ & $y^2 = x^6 + (3T + 2)x^5 + 2x^4 + 4x^3 + x^2 + 2x + 4T + 4$\\
$T + 1$ & $4T + 2$ & $4T + 2$ & $y^2 = x^6 + (3T + 4)x^5 + 2x^4 + 4x^3 + 
x^2 + 2x + 4T$\\
$T + 2$ & $4T + 1$ & $4T + 1$ & $y^2 = x^6 + (3T + 2)x^5 + 2x^4 + 4x^3 + 
x^2 + 2x + 4T + 4$\\
$T + 3$ & $4T$ & $4T$ & $y^2 = x^6 + 3Tx^5 + 2x^4 + 4x^3 + x^2 + 2x + 4T
+ 3$\\
%$T + 4$ & $2$ & $0$ & $y^2 = x^6 + 4x^5 + 4Tx^3 + 2Tx^2 + (2T + 3)x + 4T + 2$\\
$T + 4$ & $4T + 4$ & $4T + 4$ & $y^2 = x^6 + (3T + 3)x^5 + 2x^4 + 4x^3 + 
x^2 + 2x + 4T + 2$\\
$2T$ & $3T + 3$ & $3T + 3$ & $y^2 = x^6 + (T + 1)x^5 + 2x^4 + 4x^3 + x^2 +
2x + 3T + 1$\\
$2T + 1$ & $3T + 2$ & $3T + 2$ & $y^2 = x^6 + (T + 4)x^5 + 2x^4 + 4x^3 + 
x^2 + 2x + 3T$\\
$2T + 2$ & $3T + 1$ & $3T + 1$ & $y^2 = x^6 + (T + 2)x^5 + 2x^4 + 4x^3 + 
x^2 + 2x + 3T + 4$\\
$2T + 3$ & $3T$ & $3T$ & $y^2 = x^6 + Tx^5 + 2x^4 + 4x^3 + x^2 + 2x + 3T
+ 3$\\
%$2T + 4$ & $2$ & $0$ & $y^2 = x^6 + 4x^5 + 3Tx^3 + 4Tx^2 + (4T + 3)x + 3T + 2$\\
$2T + 4$ & $3T + 4$ & $3T + 4$ & $y^2 = x^6 + (T + 3)x^5 + 2x^4 + 4x^3 + 
x^2 + 2x + 3T + 2$\\
%$3T$ & $2$ & $0$ & $y^2 = x^6 + 4x^5 + (2T + 4)x^3 + (T + 2)x^2 + Tx + 2T + 1$\\
$3T$ & $2T + 3$ & $2T + 3$ & $y^2 = x^6 + (4T + 1)x^5 + 2x^4 + 4x^3 + x^2
+ 2x + 2T + 1$\\
%$3T + 1$ & $2$ & $0$ & $y^2 = x^6 + 4x^5 + (2T + 3)x^3 + (T + 4)x^2 + (T + 2)x + 2T$\\
$3T + 1$ & $2T + 2$ & $2T + 2$ & $y^2 = x^6 + (4T + 4)x^5 + 2x^4 + 4x^3 +
x^2 + 2x + 2T$\\
$3T + 2$ & $2T + 1$ & $2T + 1$ & $y^2 = x^6 + (4T + 2)x^5 + 2x^4 + 4x^3 +
x^2 + 2x + 2T + 4$\\
$3T + 3$ & $2T$ & $2T$ & $y^2 = x^6 + 4Tx^5 + 2x^4 + 4x^3 + x^2 + 2x + 
2T + 3$\\
$3T + 4$ & $2$ & $0$ & $y^2 = x^6 + 4x^5 + 2Tx^3 + Tx^2 + (T + 3)x + 2T +
2$\\
$3T + 4$ & $2T + 4$ & $2T + 4$ & $y^2 = x^6 + (4T + 3)x^5 + 2x^4 + 4x^3 +
x^2 + 2x + 2T + 2$\\
%$4T$ & $2$ & $0$ & $y^2 = x^6 + 4x^5 + (T + 4)x^3 + (3T + 2)x^2 + 3Tx + T + 1$\\
$4T$ & $T + 3$ & $T + 3$ & $y^2 = x^6 + (2T + 1)x^5 + 2x^4 + 4x^3 + x^2 + 
2x + T + 1$\\
%$4T + 1$ & $2$ & $0$ & $y^2 = x^6 + 4x^5 + (T + 3)x^3 + (3T + 4)x^2 + (3T + 2)x + T$\\
$4T + 1$ & $T + 2$ & $T + 2$ & $y^2 = x^6 + (2T + 4)x^5 + 2x^4 + 4x^3 + x^2
+ 2x + T$\\
$4T + 1$ & $4T$ & $T + 3$ & $y^2 = x^6 + 3Tx^5 + Tx^4 + Tx^3 + 2Tx^2 + 
(T + 3)x + T$\\
$4T + 2$ & $T + 1$ & $T + 1$ & $y^2 = x^6 + (2T + 2)x^5 + 2x^4 + 4x^3 + x^2
+ 2x + T + 4$\\
%$4T + 4$ & $2$ & $0$ & $y^2 = x^6 + 4x^5 + Tx^3 + 3Tx^2 + (3T + 3)x + T + 2$\\
$4T + 4$ & $T + 4$ & $T + 4$ & $y^2 = x^6 + (2T + 3)x^5 + 2x^4 + 4x^3 + x^2
+ 2x + T + 2$\\
\end{supertabular}
\end{center}
\vskip 0.5cm
\end{eg}

\begin{eg}
Let $C_f$ be the smooth hyperelliptic curve $y^2=f(x)$ over $\F_7(T)$,
where $f(x)$ is one of the polynomials in Table \ref{tbl:eg7nonss}. Then 
$$
\End(J(C_f)) = \End_{\F_7(T)}(J(C_f)) = \Z[\eta].
$$
\newpage
\begin{center}
\tablehead{
\hline
$b$ & $c$ & $d$ & $C:y^2=f_{b,c,d}(x)$ reduced modulo $7$ \\
\hline
\hline
}
\tabletail{\hline}
\topcaption{Some curves  \mbox{$C$} over  \mbox{$\F_7(T)$} with \mbox{$\End(J(C))=\Z[\eta]$}}
\begin{supertabular}
{|r@{\hspace*{0.5cm}}r@{\hspace*{0.5cm}}r@{\hspace*{0.5cm}\extracolsep{\fill}}l|}
\label{tbl:eg7nonss}
$0$ & $0$ & $T$ & $y^2 = x^6 + 2x^4 + (3T + 2)x^3 + 5x^2 + 6x + 1$\\
$0$ & $2$ & $T$ & $y^2 = x^6 + 4x^5 + 3x^4 + 3Tx^3 + 3x^2 + 3x + 1$\\
$0$ & $4$ & $T$ & $y^2 = x^6 + x^5 + 5x^4 + 3Tx^3 + 2x^2 + 1$\\
$1$ & $0$ & $5T + 4$ & $y^2 = x^6 + Tx^4 + 3Tx^3 + (T + 1)x^2 + 4x + 5$\\
$3$ & $6$ & $4T + 2$ & $y^2 = x^6 + 5x^5 + (T + 5)x^4 + Tx^2 + 5x + 6$\\
$4$ & $1$ & $3T + 1$ & $y^2 = x^6 + 2x^5 + (T + 3)x^4 + 4Tx^3 + Tx^2 + 
3$\\
$6$ & $0$ & $2T$ & $y^2 = x^6 + (T + 2)x^4 + (T + 5)x^3 + Tx^2 + x + 4$\\
$T + 3$ & $0$ & $4$ & $y^2 = x^6 + (5T + 3)x^4 + 3Tx^2 + 5Tx + 4T + 6$\\
$3T + 3$ & $3T$ & $6T + 6$ & $y^2 = x^6 + 6Tx^5 + 2Tx^4 + (2T + 4)x^2 +
5T + 6$\\
$4T + 3$ & $4T$ & $T + 6$ & $y^2 = x^6 + Tx^5 + 5Tx^4 + (5T + 4)x^2 + 2T
+ 6$\\
$5T + 3$ & $5T$ & $3T + 6$ & $y^2 = x^6 + 3Tx^5 + Tx^4 + (T + 4)x^2 + 6T
+ 6$\\
$6T + 3$ & $6T$ & $5T + 6$ & $y^2 = x^6 + 5Tx^5 + 4Tx^4 + (4T + 4)x^2 +
3T + 6$\\
\end{supertabular}
\end{center}
\end{eg}

\begin{eg}
Let $C_f$ be the smooth hyperelliptic curve $y^2=f(x)$ over $\F_{11}(T)$,
where $f(x)$ is one of the polynomials in Table \ref{tbl:eg11nonss}. Then 
$$
\End(J(C_f)) = \End_{\F_{11}(T)}(J(C_f)) = \Z[\eta].
$$
\begin{center}
\tablehead{
\hline
$b$ & $c$ & $d$ & $C:y^2=f_{b,c,d}(x)$ reduced modulo $11$ \\
\hline
\hline
}
\tabletail{\hline}
\topcaption{Some curves  \mbox{$C$} over  \mbox{$\F_{11}(T)$} with \mbox{$\End(J(C))=\Z[\eta]$}}
\begin{supertabular}
{|r@{\hspace*{0.5cm}}r@{\hspace*{0.5cm}}r@{\hspace*{0.5cm}\extracolsep{\fill}}l|}
\label{tbl:eg11nonss}
$0$ & $0$ & $T + 6$ & $y^2 = x^6 + 2x^4 + 7Tx^3 + 5x^2 + 6x + 1$\\
$0$ & $8$ & $T + 9$ & $y^2 = x^6 + 5x^5 + 5x^4 + 7Tx^3 + 2x^2 + 1$\\
$5$ & $10$ & $6T + 2$ & $y^2 = x^6 + 9x^5 + (T + 5)x^4 + Tx^2 + 9x + 10$\\
$T + 5$ & $0$ & $6$ & $y^2 = x^6 + (9T + 3)x^4 + 10Tx^2 + Tx + 4T + 10$\\
$2T + 5$ & $10T$ & $7T + 10$ & $y^2 = x^6 + 9Tx^5 + 9Tx^4 + (9T + 8)x^2 + 8T + 10$\\
$3T + 1$ & $2T$ & $4T + 6$ & $y^2 = x^6 + 4Tx^5 + 4Tx^4 + 4x^2 + (7T + 7)x + T + 5$\\
%$3T + 5$ & $0$ & $0$ & $y^2 = x^6 + 2x^4 + Tx^3 + (3T + 10)x^2 + 3Tx + T + 10$\\
$3T + 10$ & $4T$ & $5T + 6$ & $y^2 = x^6 + 8Tx^5 + 4x^4 + 6Tx^3 + 6x^2 
+ 5x + T + 8$\\
$4T + 5$ & $9T$ & $3T + 10$ & $y^2 = x^6 + 7Tx^5 + 7Tx^4 + (7T + 8)x^2 
+ 5T + 10$\\
\end{supertabular}
\end{center}
\end{eg}

%=========================================================================
\section*{Acknowledgments}
%=========================================================================
I would like to thank my advisor, Yuri Zarhin, for his guidance.

%=========================================================================
%\section*{List of Notation}
%=========================================================================
%\noindent
%\begin{tabular}{ll}
%$D$			& real quadratic number field\\
%$d$ 			& the reduced discriminant of $D$\\
%$\omega$ 	& $\sqrt{d}$\\
%$\O$ 		& an order of $D$\\
%$\eta$		& generator of $\O$ different from $1$\\
%$c$			& the conductor of $\O$\\
%$\discr(\ldots)$	& the discriminant of\\
%$C_f$		& hyperelliptic curve $y^2=f(x)$\\
%$X=J(C)$ 		& jacobian variety of a curve $C$\\
%$\i$ 			& map $\O \hookrightarrow \End_K(X)$
%			   (extends to $D \hookrightarrow \End_K^0(X)$)\\
%$\End^0(X,\i)$	& centralizer of $\i(D)$ in $\End^0(X,\i)$\\
%$T_2(X)$		& $2$-adic Tate module of $X$\\
%$V_2(X)$		& $\Q_2$-Tate module of $X$\\
%\end{tabular}
%=========================================================================


\begin{thebibliography}{99}

\bibitem{MAGMA}
W.~Bosma, J.~Cannon,  C.~Playoust,
{\em The MAGMA algebra system, I: The user language}. 
J. Symb. Comp., {\bf 24} (1997),
235-265. (See also the MAGMA homepage at \url{http://www.maths.usyd.edu.au:8000/u/magma/})

\bibitem{BruRankJ} A.~Brumer,
{\em The Rank of $J_0(N)$}.
Ast{\'e}risque {\bf 228} (1995), 
41--68.

\bibitem{CuMRT} C.~W.~Curtis, I.~Reiner,
Methods of Representation Theory, Volume 1,
John Wiley and Sons, New York,
1981.

\bibitem{DixPG} J.~D.~Dixon, B.~Mortimer,
Permutation Groups, 
Graduate Texts in Mathematics {\bf 163},
Springer-Verlag, Berlin-Heidelberg-New York, 
1996.

\bibitem{FeitSI} W.~Feit,
{\em The computation of some Schur indices}.
Israel J. Math. {\bf 46} (1983), % No. 4, 
274--300.

\bibitem{GorFSG} D.~Gorenstein,
Finite Simple Groups, An Introduction to their Classification,
Plenum Press, New York and London,
1982.

\bibitem{HashiOnBru} K.~Hashimoto,
{\em On Brumer's family of RM-curves of genus two}.
Tohoku Math. J. {\bf 52} (2000), 
475--488.

\bibitem{IbuSC} T.~Ibukiyama, T.~Katsura, F.~ Oort,
{\em Supersingular curves of genus two and class numbers}.
Compositio Math. {\bf 57} (1986), 
127--152.

\bibitem{ManCFG} Yu.~I.~Manin, 
{\em Theory of commutative formal groups over fields of finite characteristic}. 
Uspehi Mat. Nauk, {\bf 18} (1963),
3--90

\bibitem{MasserSHJ}  D.~Masser,
{\em Specialization of some hyperelliptic jacobians}.
In: Number Theory in Progress (eds. K.~Gy{\"o}ry, H.~Iwaniec, J.~Urbanowicz), 
vol. I, pp. 293--307; de Gruyter, Berlin-New York,
1999.

\bibitem{MoriER1}  Sh.~Mori, 
{\em The endomorphism rings of some abelian varieties}.
Japanese J. Math. {\bf 2} (1976), 
109--130.

\bibitem{MoriER2}  Sh.~Mori, 
{\em The endomorphism rings of some abelian varieties.II}.
Japanese J. Math. {\bf 3} (1977),
105--109.

\bibitem{MortMPR} B.~Mortimer, 
{\em The modular permutation representations of the known doubly transitive groups}.
Proc. London Math. Soc. (3) {\bf 41}  (1980), 
1--20.

\bibitem{MumAV} D.~Mumford,
Abelian varieties, Second edition, 
Oxford University Press, London,
1974.

\bibitem{NygCC} N.~O.~Nygaard, 
{\em Slopes of powers of Frobenius on crystalline cohomology}.
Ann. Sci. Ecole Norm. Sup. {\bf 14} (1981), 
369--401.

\bibitem{OortEnd} F.~Oort,
{\em Endomorphism Algebras of Abelian Varieties}.
Algebraic Geometry and Commutative Algebra (1987),
469--502.

\bibitem{PARI2}
PARI/GP, version {\tt 2.1.5}, Bordeaux, 2004, \url{http://pari.math.u-bordeaux.fr/}.

\bibitem{SerTateRed} J.-P.~Serre, J.~Tate, 
{\em Good reduction of abelian varieties},
Annals of Mathematics {\bf 88} (1968), 
492--517.

\bibitem{SilMS} A.~Silverberg, Yu.~G.~Zarhin, 
{\em Variations on a theme of Minkowski and Serre},
J.\ Pure and Applied Algebra {\bf 111}  (1996),  
285--302.

\bibitem{SilEllCurves} J.~Silverman,
The Arithmetic of Elliptic Curves,
Graduate Texts in Mathematics {\bf 106},
Springer-Verlag, New York-Berlin-Heidelberg-Tokyo, 
1986.

\bibitem{StPM} W.~Stein, 
PARI \& MAGMA Calculator, \url{http://modular.fas.harvard.edu/calc/}.

\bibitem{VigQA} M.-F.~Vign{\'e}ras, 
Arithm{\'e}tique des Alg{\`e}bres de Quaternions, 
Lecture Notes in Mathematics {\bf 800},
Springer-Verlag, Berlin-Heidelberg-New York,
1980.

\bibitem{WilEM} J.~Wilson, 
{\em Explicit moduli for curves of genus $2$ with real multiplication by $\Q(\sqrt{5})$},
Acta. Arith. {\bf 93} (2000), 
121--138.

\bibitem{WilRM} J.~Wilson,
{\em Curves of genus $2$ with real multiplication by
a square root of $5$},
Oxford University, Ph.D. thesis,
1998.

\bibitem{YuiJV} N. ~Yui, 
{\em On the jacobian varieties of hyperelliptic curves over 
fields of characteristic $p>2$}.
Journ. of Alg. {\bf 52} (1978), 
378--410.

\bibitem{ZarCM} Yu.~G.~Zarhin,
{\em Hyperelliptic jacobians without complex multiplication}.
Math. Res. Letters {\bf 7} (2000), 
123--132.

\bibitem{ZarMR} Yu.~G.~Zarhin,
{\em Hyperelliptic jacobians and modular representations}.
In: Moduli of abelian varieties,
(C. Faber, G. van der Geer, F. Oort, eds.), pp. 473--490, 
Progress in Math., Vol. {\bf 195},
Birkhauser, Basel-Boston-Berlin, 
2001.

\bibitem{ZarVS2R} Yu.~G.~Zarhin,
{\em Very simple $2$-adic representations and hyperelliptic jacobians}.
Moscow Math. J. {\bf 2} (2002), issue 2,
 403--431.

\bibitem{ZarCMP} Yu.~G.~Zarhin,
{\em Hyperelliptic jacobians without 
complex multiplication in positive characteristic}, 
Math. Res. Letters {\bf 8} (2001), 
429--435.

\bibitem{ZarNS} Yu.~G.~Zarhin,
{\em Non-supersingular hyperelliptic jacobians},
Bull. Soc. Math. France, {\bf 132} (2004)
617--634.

\end{thebibliography}
\end{document}